\input eplain 
\magnification  1200
\ifx\eplain\undefined \input eplain \fi

\baselineskip13pt 
\vsize=9.1truein

\overfullrule=0pt
\def\zero{{\bf 0}}
 

\def\wtilde{\widetilde} 
\font\mbi=msbm7
\def\sZ{{\hbox{\mbi\char 90}}}

\font\scaps=cmcsc10

\font\smalsmalbf=cmbx8
\font\tenbi=cmmib10

\font\eightbi=cmmib9

\font\fivebi=cmmib5
\newfam\bmifam\textfont\bmifam=\tenbi\scriptfont
\bmifam=\eightbi\scriptscriptfont\bmifam=\fivebi

\font\smalltenrm=cmr8
\font\smallteni=cmmi8
\font\smalltensy=cmsy8
\font\smallsevrm=cmr6   \font\smallfivrm=cmr5
\font\smallsevi=cmmi6   \font\smallfivi=cmmi5
\font\smallsevsy=cmsy6  \font\smallfivsy=cmsy5
\font\smallsl=cmsl8      \font\smallit=cmti8

\def\smallfonts{\lineadj{80}\textfont0=\smalltenrm  \scriptfont0=\smallsevrm
                \scriptscriptfont0=\smallfivrm
    \textfont1=\smallteni  \scriptfont1=\smallsevi
                \scriptscriptfont0=\smallfivi
     \textfont2=\smalltensy  \scriptfont2=\smallsevsy
                \scriptscriptfont2=\smallfivsy
      \let\it\smallit\let\sl\smallsl\smalltenrm}

\font\eightbi=cmmib10 at 8pt

\font\smathbold=msbm8\font\ssmathbold=msbm5
\font\mathbold=msbm9 at 10pt\def\Hn{{\hbox{\mathbold\char72}}^n}
\def\sHn{{{\hbox{\smathbold\char72}}_{}^n}}\def\ssHn{{{\hbox{\ssmathbold\char72}}_{}^n}}
\def \C{{\hbox{\mathbold\char67}}}
\def\N{{\hbox{\mathbold\char78}}}
\def\R{{\hbox{\mathbold\char82}}}\def\Z{{\hbox{\mathbold\char90}}}
\def\sR{{\hbox{\smathbold\char82}}}

\def\e{\epsilon}

\def\d{\delta}


\def\imp{{\;\;\Longrightarrow\;\;}}
\def\lineadj#1{\normalbaselines\multiply\lineskip#1\divide\lineskip100
\multiply\baselineskip#1\divide\baselineskip100
\multiply\lineskiplimit#1\divide\lineskiplimit100}

\def\remark#1.{\medskip{\noin\bf Remark #1.\enspace}}
\def\endpf{$$\eqno/\!/\!/$$}

\def\pf#1.{\smallskip\noin{\bf  #1.\enspace}}

\def\noin{\noindent}

\def\esup{\sup_{x\in N}}

\def\ds{\displaystyle}
\def\ts{\textstyle}

\def\e{\epsilon}\def\part{\partial_t}

\def\isn{\int_{S^{n-1}}}\def\wtilde{\widetilde}

\def\Rn{\R^n}
\def\im{\int_M}
\def\irn{\int_{\sR^n}}
\def\ref#1{{\bf{[#1]}}}

\def\essup{{\rm {ess\hskip.2em sup}}}
\def\p{\partial}

\def\essinf{{\rm {ess\hskip2pt inf }}\;}

\def\half{{{1\over2}}}
\def\A{{\bf A}}\def\G{{\cal G}}

\def\s{\sigma}
\def\supp{{\rm supp }\,}
\def\H{{\hbox{\mathbold\char72}}}\def\b{\beta}
\def\a{\alpha}
\def\na{{n/\a}}\def\nsa{{n\over \a}}
\def\nna{{n\over n-\a}}
\def\bp{{\b'}}\bigskip\def\b{\beta}\def\bi{{1\over\b}}\def\bpi{{1\over\b'}}


\centerline{\bf Adams inequalities for Riesz subcritical potentials}

\bigskip \centerline{Luigi Fontana, Carlo Morpurgo}
\footnote{}{\smallfonts
\hskip-2.4em  The work of the second author was partially supported by NSF Grant DMS-1401035 and by Simons Foundation Collaboration Grant 279735}
\midinsert
{\smalsmalbf Abstract. }{\smallfonts 
We derive Adams inequalities for potentials on  general measure spaces, extending and improving previous results in  [FM1]. The integral  operators involved, which we call ``Riesz subcritical", have  kernels whose decreasing rearrangements are not worse than that of the Riesz kernel  on $\sR^n$, where the kernel is large, but they behave  better where the kernel is small.  The new element is a ``critical integrability" condition on the kernel at infinity.  Typical examples of such kernels are fundamental solutions of  nonhomogeneous differential, or pseudo-differential, operators. Another example is the Riesz kernel itself restricted to suitable measurable sets, which we name ``Riesz subcritical domains". Such domains are characterized in terms of their  growth at infinity.  As a consequence of the general results we obtain several new sharp Adams and Moser-Trudinger inequalities on $\sR^n$, on the hyperbolic space,  on Riesz subcritical domains, and on domains where the Poincar\'e inequality holds.
}

\endinsert\medskip\bigskip
\def\am{{\alpha\over2}}
\centerline{\scaps 1. Introduction}\bigskip
 It is well understood by now that the validity of the  Moser-Trudinger inequality relative to a sufficiently well-behaved differential, or pseudodifferential, operator $P$ is strongly related to the behavior of the fundamental solution of $P$ around its singularity. More  specifically, if $E$ and $F$ are open sets with {\it finite} Lebesegue measure on $\Rn$ and $P$ is of order $\alpha\in (0,n)$  and invertible, then one has the Moser-Trudinger inequality 
$$\int_E \exp\bigg[\gamma|u(x)|^\nna\bigg]dx\le  C\qquad u\in W_0^{\a,\nsa}(F),\qquad \|Pu\|_\na\le 1\eqdef{mt1}$$
for some $C$ depending  on $|E|,|F|, n,\a$, where the best constant $\gamma=\gamma(P)$ is related to the leading term of the kernel of $P^{-1}$ around the diagonal, and ultimately to the homogeneous principal symbol of $P$. This has been analyzed in great generality in our earlier paper [FM1], but the first important result goes back to Adams, who considered the case $P=\nabla^\alpha$, $\alpha$ an integer, and $E=F$. Here $\nabla^\a=(-\Delta)^{\am}$ for $\alpha$ even and $\nabla^\a=\nabla(-\Delta)^{\alpha-1\over2}$ if $\alpha$ odd. When $\alpha$ is even then the inverse of $\nabla^\a$ is of course the  Riesz potential  $c_\a I_\a f$ where
$$I_\a f=|x|^{\a-n}*f,\qquad 
c_\a={\Gamma\big({n- \a\over2}\big)\over 2^\a\pi^{n/2}\Gamma\big({ \a\over2}\big)}.\eqdef{calpha}$$
In [A1] Adams proved the following  sharp inequality (which we  call ``Adams inequality")
$$\int_E \exp\bigg[{1\over |B_1|} |I_\a f(x)|^\nna\bigg]dx\le  C|E|\qquad f\in L^\na(E),\qquad \|f\|_\na\le 1\eqdef{mt2}$$
where $B_1$ is the unit ball of $\Rn$. From \eqref{mt2} Adams  derived the sharp form of \eqref{mt1},   when $E=F$, for the operators $P=\nabla^\a$, with sharp exponential constant $\gamma$ given as  
$$\gamma(\nabla^\a)=\cases{\ds{c_\a^{-\nna}\over|B_1|} & if $\a$ even\cr \ds{\big((n-\a-1)c_{\a+1}\big)^{-\nna}\over|B_1|}& if $\a$ odd.\cr}$$
In [FM1] we pushed Adams' argument to the extreme, by considering general integral operators 
$$Tf(x)=\int_M k(x,y)f(y) d\mu(y),\qquad x\in N\eqdef T$$ on spaces of finite measure $(M,\mu),\;(N,\nu)$, and by showing that  the validity the estimate
$$\int_E \exp\bigg[\gamma |T f(x)|^\nna\bigg]dx\le  C\qquad f\in L^\na(F),\qquad \|f\|_\na\le 1\eqdef{mt3}$$
is reduced 
 to a couple of growth estimate on the kernel $k$. To describe such estimates we introduce the partial nonincreasing rearrangements 
 
 $$k_1^*(x,t)=[k(x,\cdot)]^*(t),\qquad k_2^*(y,t)=k[(\cdot,y)]^*(t)$$ 
 where $k[(x,\cdot)]^*$ is the nonincreasing rearrangement of $k(x,y)$ with respect to $y$ on $(M,\mu)$ for fixed $x\in N $, and $k[(\cdot,y)]^*$ is its analogue on $(N,\nu)$ for fixed $y\in M$.  The corresponding  maximal nonincreasing rearrangements are defined as 
 $$k_1^*(t):=\sup_{x\in N} k_1^*(x,t),\qquad
k_2^*(t):=\sup_{y\in M} k_2^*(y,t).$$
where (here and throughout the paper)  $``\sup"$ means ``\essup" and $\inf$ means ``$\essinf$".\smallskip

 The main Theorem in [FM1] states that if for some $A,B>0$, $\gamma>1$ and $\sigma\in (0,1]$ there holds 
 $$k_1^*(t)=t^{-{n-\a\over n}}(A^{n-\a\over n}+O(|\log t|)^{-\gamma}\big),\qquad 0<t\le1\eqdef{106a}$$
$$k_2^*(t)\le Bt^{-{n-\a\over\sigma n}},\qquad t>0\eqdef{107a}$$
 then \eqref{mt3} with $E=N,F=M$ holds with $\gamma=\sigma A^{-1}$, which is generally best possible. The growth estimates \eqref{106a},\eqref{107a} reflect the nature of the singularity of fundamental solutions of invertible elliptic (pseudo) differential operators, the Riesz potential being the model case. The constant $\sigma$ in \eqref{107a} allows from some flexibility in the choice of measure in the target space $(N,\nu)$, on $\Rn$ typically singular or Hausdorff measures.
 
 As an application we obtained  \eqref{mt1} 
for a general class of smooth invertible  elliptic operators $P$ on bounded domains $E=F=\Omega$. In such setting one has that $P^{-1}f$ can be written as an integral operator with a kernel $k$ satisfying  $k(x,y)\sim g_\a (x,x-y)$, when $y$ is close to $x$, and where $g_\a$ is homogeneous of order $\a-n$ in the second variable. This local asymptotic behavior is all that is needed in order to obtain  \eqref{mt1} with (generally sharp)   exponential constant
$$\gamma=\bigg({1\over n}\sup_{x\in \Omega}\int_{S^{n-1}}\bigg|\Big({1\over p_\a(x,\cdot)}\Big)^\wedge(\omega)\bigg|^\nna d\omega\bigg)^{-1}\eqdef {gamma}$$
where $p_\a(x,\xi)$ is the strictly homogeneous principal symbol of $P$ (see [FM1, Thm. 10]).
 
Things are more complicated if we ask for the constant $C$ in either \eqref{mt1} or \eqref{mt3} to be independent of the measure of $F$; this means for example that in \eqref{mt1} we are only requiring $u\in W^{\a,\nsa}(\Rn)$. For operators $P$ which are homogeneous of order $\alpha$, or for convolution kernels $K$ which are homogeneous of order $\a-n$ there is no hope of obtaining uniformity on $|F|$ in the above inequalities for any value of $\gamma$ - this is seen by a simple dilation argument. In these cases one can overcome the lack of control on the measure of the support by imposing  norm additional conditions on $\|u\|_\na$ for \eqref{mt1},  or  on $\|Tf\|_\na$ for \eqref{mt3}. In the case  $P=\nabla^\a$ any of the norms
$$\|u\|_{\na,q}:=\cases{\big(\|u\|_\na^{q\na}+\|\nabla^\a u\|_\na^{q\na}\big)^{\a\over qn} & if $1\le q<\infty$\cr \max\big\{\|u\|_\na,\|\nabla^\a u\|_\na\big\} & if $ q=+\infty$\cr}\eqdef{norm}$$ 
is equivalent to the Sobolev norm in $ W^{\a,\nsa}(\Rn)$.   Under the condition $\|u\|_{\na,q}\le 1$, inequality \eqref{mt1} holds with sharp constant $\gamma=\gamma(\nabla^\a)$ when $q= 1$ (see [FM2], [LL], [LR], [MaS], [R]), (in [MaS]  an even stronger Moser-Trudinger inequality is proven) and for any $\gamma<\gamma(\nabla^\a)$ when $q>1$ (see [AT], [C], [do\'O], [FM2], [P]). When  $P=(-\Delta)^{\a\over2}$, any real $\a\in (0,n)$, and other homogeneous $P$  (see [FM2, Theorem 3]),
and for Riesz-like potentials in the Adams inequality \eqref{mt3},  under the condition $\|f\|_\na^\na+\|Tf\|_\na^\na\le 1$ ([FM2, Thm 5]).

\smallskip

As it turns out the above-mentioned difficulties associated with homogeneous operators are due to insufficient decay at infinity of their fundamental solutions. In the case of  the Bessel operator $P=(I-\Delta)^{\a/2}$, for example, the result in \eqref{mt1} holds for $F=\Rn$, and $E$ with finite measure, with the same  sharp constant  $\gamma(\nabla^\a)$. This was proved first by Adams himself in his original paper [A, Thm. 3] in the case $\alpha=2$, by writing $u$ as a Bessel potential $u=G_\a*f$, with $f=(I-\Delta)^{\a/2}u$, and by proving  the sharp inequality \eqref{mt2}, with $Tf=G_\a*f$. Such a proof was just a small modification of the proof he gave for the Riesz potential, which was possible given the good exponential decay of $G_\a$ at infinity. We should mention here that there was nothing peculiar about $\a=2$ in that argument, except the special form taken by $\|(I-\Delta)u\|_2$, and that the proof could have been carried out in full generality (this has been done in [LL2], and a different proof for $\a$ integer appears in [RS].)

Upon reading Adams' result we realized that the entire measure-theoretic machinery developed in our paper [FM1] could be extended to incorporate inequalities such as \eqref{mt3} in the case of $F$ having possibly infinite measure. We achieved this by adding to estimates \eqref{106a}, \eqref{107a}  the following  {\it critical integrability} condition  on the kernel $k$:
$$\sup_{x\in N} \int_1^\infty k_1^*(x,t)^\nna dt<\infty.\eqdef{108a}$$
In this paper we show that condition \eqref{108a} is sufficient and  also essentially necessary in order for \eqref{mt2} to hold (see Theorem 1 and Theorem 6).  

The proof of this result is based on several improvements of the arguments given in [FM1], including a new improved version of  O'Neil's lemma.

We will call a kernel $k(x,y)$, and its corresponding potential, {\it Riesz subcritical} if it satisfies the integrability condition \eqref{108a} together with 
the estimates $$k_1^*(t)\le C t^{-{n-\a\over n}},\qquad k_2^*(t)\le Bt^{-{n-\a\over \sigma n}},\qquad t>0\eqdef{conditions}$$
for some fixed constants $C,B$, and for some $\sigma\in (0,1]$.
Likewise,  a kernel (and its potential) will be called {\it Riesz critical} if it satisfies \eqref{conditions} but   the critical integrability  condition \eqref{108a} does not hold.

     For example, when $M=\Rn$ and $N\subseteq\R^n$ with the Lebesgue measure, the Riesz potential itself is Riesz critical. Still on $\R^n$, if $T=K*f$, then $K$ satisfies \eqref{106a} and \eqref {107a}, hence \eqref{conditions}, with $\sigma=1$,  if it behaves like a Riesz kernel near 0. On the other hand,  $K$ satisfies the critical integrability condition \eqref{108a} if, loosely speaking,
it decays at infinity ever so slightly better than the Riesz  kernel, in the sense that   $K$ is in $L^{\nna}$ and bounded outside a large ball. For such potentials the Adams inequality holds, on $\R^n$,  in the same spirit as the original Adams result (see Theorem 7).  
The main examples of kernels of this type are those arising as fundamental solutions of invertible elliptic {\it non-homogeneous} differential operators. We have already mentioned above that invertible elliptic differential operators have kernels which are locally behaving like a Riesz kernel, and if they are homogeneous of order $\alpha$ then they are Riesz critical. However, in section 4 we will show that if $P$ is in a large class of non homogeneous elliptic operators with  constant coefficients, then its fundamental solution is indeed Riesz subcritical. For such operators we then have a Moser-Trudinger inequality  of type \eqref{mt1} which hold for $u\in W_0^{\a,\nsa}(\Rn)$ and under the condition $\|Pu\|_\na\le1$ (Theorem 14).

Other examples of differential operators with Riesz subcritical inverses are the powers of the Laplacian in hyperbolic space. In Theorem 15 we will give a sharp version of the Moser-Trudinger inequality for such operators, extending the known results in the case $\alpha=1$ to higher powers.
The techniques developed in this paper could very likely  be applied  to other noncompact manifolds. 

An interesting question  arising from this work is the  following: on which measurable sets  $\Omega$ is the Riesz potential itself also  Riesz subcritical? This is the same as asking Riesz subcriticality of the kernel $\chi_\Omega^{} (x)|x-y|^{\a-n}\chi_\Omega^{}(y)$. We  call such sets {\it Riesz subcritical domains}. It turns out that the condition for subcriticality is quite explicit, and also independent of $\alpha$: if $\Lambda_\Omega(r)= |\Omega\cap B(x,r)|$, then $\Omega$ is Riesz subcritical if and only if  
$$\sup_{x\in \Omega} \int_1^\infty{\Lambda_\Omega(x,r)\over r^{n+1}}dr<\infty. \eqdef{RSC}$$

In essence, $\Omega$ is  Riesz subcritical if it ``loses enough dimensions at infinity''. Examples include sets of finite measure, and  sets which are bounded in one or more dimensions (like ``strips'').

On such subcritical sets $\Omega$ the Adams inequality \eqref{mt2} holds for all $f\in L^\na(\Omega)$ under the sole condition that $\|f\|_\na\le1$, as in the original Adams result (see Theorem 9). In some sense Riesz subcritical domains are the best replacements of sets with finite measure, as far as the Adams inequality is concerned.  In section 3 we will give more examples of  Riesz subcritical sets, and we will also show that \eqref{RSC} is  a necessary condition for the Adams inequality to hold (see Theorem 9).

In section 6 we explore some relations between the Moser-Trudinger and the Poincar\'e inequalities. This connection comes about from the so-called {\it regularized  form} of the Moser-Trudinger inequality for $\nabla^\a$, on sets of infinite measure, which takes the form 
$$\int_{\Omega}\exp_{[\nsa-2]}\Big[\gamma(\nabla^\a)|u(x)|^\nna\Big]dx\le C(1+\|u\|_{L^\na(\Omega)}^\na),\qquad \|\nabla^\a u\|_{L^\na(\Omega)}\le 1,\eqdef{regmt}$$
and which we prove to be  valid when $\Omega$ is an open Riesz subcritical set in $\R^n$ (see Corollary~10).  If in addition  the Poincar\'e inequality 
$\|u\|_\na\le C\|\nabla^\a u\|_\na$ is valid on $W_0^{\a,\na}(\Omega)$, then the exponential integral in \eqref{regmt} is clearly uniformly bounded. The question is whether or not such uniformity of the exponential integral holds in domains where the  above Poincar\'e inequality holds, but which are not Riesz subcritical. Examples of such domains are complements of periodic nets of closed balls of fixed radius.  In [BM], Battaglia e Mancini proved that indeed this is the case when $\alpha=1$.  While we were not able to prove that the same is true for all integer $\alpha$, we instead show that on a large class of Poincar\'e domains (not necessarily subcritical) the Moser-Trudinger inequality does hold. Such domains are characterized by a geometric  condition given in terms of  the {\it strict inradius}, term introduced by P. Souplet in [So], based on ideas due to Agmon [Ag].  The proof of our result is a nice application of the general measure-theoretic Adams inequality. 

\bigskip

\centerline{\scaps 2. Adams inequalities for Riesz subcritical potentials}\bigskip

 Suppose that  $(M,\mu)$ and $(N,\nu)$ are measure spaces  and  that  $T$   is an integral operator of type
$$Tf(x)=\int_M k(x,y)f(y)d\mu(y),\qquad x\in N,\eqdef T$$
with $k:N\times M\to[-\infty,\infty]$ a $\nu\times\mu$-measurable function. Define the partial distribution functions of $k$ as 
 $$\cases{\lambda_1(x,s)=\mu(\{y\in M:\;|k(x,y)|>s\}) & $\;x\in N$\cr \cr
 \lambda_2(y,s)=\nu(\{x\in N:\;|k(x,y)|>s\}) &  $\;y\in M.$\cr}\eqdef{pdis1}.$$
 The corresponding partial nonincreasing  rearrangements are given as 
$$\cases{k_1^*(x,t)=\inf\{s>0:\, \lambda_1(x,s)\le t\} & $\; x\in N$ \cr\cr
k_2^*(y,t)=\inf\{s>0: \, \lambda_2(y,s)\le t\} & $\; y\in M$\cr} \eqdef{pdis2}$$
and the maximal nonincreasing rearrangements are defined as
$$k_1^*(t)=\sup_{x\in N} k_1^*(x,t),\qquad
k_2^*(t)=\sup_{y\in M} k_2^*(y,t)\eqdef{pdis3}$$
where, once again, sup and inf are in the sense of essential sup and essential inf.

\proclaim Theorem 1. Suppose that there exist constants $A,B,H>0$ such that for some $\tau>0$ and $\gamma>1$
$$k_1^*(t)\le A^{1\over\beta}t^{-{1\over\beta}}\big(1+H(1+|\log t|)^{-\gamma}\big),\qquad 0<t\le\tau\eqdef{106}$$
$$k_2^*(t)\le Bt^{-{1\over\sigma \beta}},\qquad t>0\eqdef{107}$$
$$J_\tau:={1\over A}\,\sup_{x\in N}\, \int_{\tau}^\infty\big (k_1^*(x,t)\big)^{\b}dt<\infty,\eqdef{108}$$
for some $\b>1$ and $0<\sigma\le 1$.
Then,  $Tf$ is well defined and  finite a.e. for $f\in L^\bp(M)$, and  there exists a constant $C=C(\b,\sigma,\gamma,A,B,H)$ such that for each $f\in L^{\b'}(M)$ with $\|f\|_{\b'}\le 1$
and for each measurable $E\subseteq N$ with $\nu(E)<\infty$ we have
$$\int_E \exp\bigg[{\sigma\over A}|Tf(x)|^\beta\bigg]d\nu(x)\le  Ce^{\sigma J_\tau}\big(\nu(E) +\tau^\sigma(1+J_\tau)\big).\eqdef{109}$$
Moreover,   for all $f\in L^\bp(M)$ such that $\|f\|_\bp\le1$  we have  
$$\int_N \exp_{[\b'-2]}\bigg[{\sigma\over A}|Tf(x)|^\beta\bigg]d\nu(x)\le  Ce^{\sigma J_\tau}\big(\|Tf\|_\bp^\bp+\tau^\sigma(1+J_\tau)\big),\eqdef{110}$$
where  the  ``regularized exponential'' function $\exp_m$ is defined as
$$\exp_m(t)=e^t-\sum_{k=0}^m{t^k\over k!},\qquad m=0,1,...$$
and where $[\lambda]$ denotes the smallest integer  greater or equal to $\lambda\in \R.$
\par

Later in Theorem 6 we will show that for reasonable kernels the critical integrability  condition \eqref{108} is also necessary for the Adams inequality to hold.  We note that the earlier version of this theorem appearing  in [FM3, Thm. 3], has a slightly stronger assumption, namely 
$$\int_1^\infty \big(k_1^*(t)\big)^\b dt<\infty.\eqdef{stronger}$$
Assumption \eqref{stronger} is in general easier to verify than \eqref{108}. Moreover Theorem 1 can be  deduced from \eqref{stronger} by  using 
a milder modification of O'Neil's Lemma than the one   we need (see Lemma 5 below)  if we only assume \eqref{108} (see  Lemma 2 in [FM1] and Note 3 after Lemma 5). 
On the other hand, there are examples for which \eqref{108} holds  but \eqref{stronger} does not.  \eject
It is not too hard to produce kernels for which

$$\sup_{x\in N}\int_1^\infty \big(k_1^*(x,t)\big)^\b dt<\infty,\qquad \int_1^\infty \big(k_1^*(t)\big)^\b dt=+\infty.\eqdef{supint3}$$

For a somewhat artificial kernel, consider on $\R^n$ with the Lebesgue measure 
$$k(x,y)=\cases{|y|^{\a-n} & if $\ds{|x|\over2}\le |y|<|x|$\cr 
0 & otherwise\cr}$$
and let $Tf(x)=\int_{\sR^n} k(x,y)f(y)dy.$ In this case we have
$$k_1^*(x,t)=\cases{\bigg(\ds{t\over |B_1|}+\Big(\ds{|x|\over2}\Big)^n\bigg)^{-{n-\a\over n}} & if $0<t< |B_1| |x|^n(1-2^{-n})$\cr 0 & if $t\ge |B_1||x|^n(1-2^{-n})$\cr}\le |B_1|^{n-\a\over n}t^{-{n-\a\over n}} =k_1^*(t)$$
which implies \eqref{supint3}, with $\b={\nna}$, and it is easy to check that $k_2^*(t)\le Ct^{-{n-\a\over n}}$.  The Adams inequality trivially holds, since by H\"older's inequality $|Tf(x)|\le C\|f\|_\na$  for all $x$. Note in passing that the same example shows that $T$ is not continuous from $L^{\na}(\Rn)$ to any $L^p(\Rn)$ (consider the family $f_r(y)=|y|^{-\a}\chi_{r/2\le|y|<r}^{}(y))$.

A more interesting example satisfying \eqref{supint3} is discusssed  in \eqref{dominio1}, where $k(x,y)=\break|x-y|^{\a-n}$ and 
$M=N= \Omega$ for a suitable $\Omega$ of infinite measure, whose construction is not trivial.
Adams' inequality for $\Omega$ follows immediately from Theorem 1  in the present form,
whereas the previous version in [FM3] based on assumption \eqref{stronger} is unable to give 
a sensible answer for $\Omega$.

 The right hand sides of \eqref{109}, \eqref{110} are also improvements of previous versions appearing in [FM3, Theorem 3].
\smallskip

If we constrain our functions $f$ to be supported on a given set $F\subseteq M$ with finite measure   we obtain the following refinement of Theorem 1 in [FM1]:

\proclaim Corollary 2. Suppose that 
 $$k_1^*(t)\le A^{1\over\beta}t^{-{1\over\beta}}\big(1+H(1+|\log t|)^{-\gamma}\big),\qquad \forall t>0\eqdef{106-a}$$
$$k_2^*(t)\le Bt^{-{1\over\sigma \beta}},\qquad \forall t>0.\eqdef{107-b}$$
Then,  there exists a constant $C=C(\b,\sigma,\gamma,A,B,H)$ such that for each measurable $F\subseteq M$ and each  $f\in L^{\b'}(F)$ with $\|f\|_{L^{\b'}(F)}\le 1$ and for each measurable $E\subseteq N$ 
$$\int_E \exp\bigg[{\sigma\over A}|Tf(x)|^\beta\bigg]d\nu(x)\le  C\big(\nu(E) +\mu(F)^\sigma\big).\eqdef{109-a}$$\par

\smallskip
\noin{\bf Note.} When $F\subseteq M$ is measurable,  the space $L^p(F)$ is defined as the space of those measurable $f:E\to \R$ such that their zero extension to $M$ is in $L^p(M)$. Equivalently, it is the space of functions $f\in L^p(M)$ which are 0 a.e. outside $F$.
\smallskip
\pf Proof  of Corollary 2.  If $\mu(F)<\infty$ then we can apply Theorem 1 to the measurable space $(F,\mu)$ as a subspace of  $M$. If $k_{1,F}^*(x,t),\, k_{1,F}^*(t)$ denote the rearrangements with respect to $F$ then 
$$k_{1,F}^*(t)\le k_1^*(t)\le A^{1\over\beta}t^{-{1\over\beta}}\big(1+H(1+|\log t|)^{-\gamma}\big),\qquad 0<t\le \mu(F)$$
(with $H$ independent of $F$!) and $k_{1,F}^*(t)=0$ for $t>\mu(F)$. Hence condition \eqref{108} is verified for $\tau=\mu(F)$ which gives the right hand side in \eqref{109-a}.\endpf

When $T$ is the Riesz potential Corollary 2 gives the following refinement of Adams' inequality \eqref{mt2} and of Theorem 7 in [FM1]:

\proclaim Corollary 3.  Let  $\nu$ be a positive Borel measure on $\Rn$ satisfying
 $$\nu\big(B(x,r)\big)\le Q r^{\sigma n},\qquad \forall x\in \Rn,\,\forall r>0\eqdef{b}$$
for some   $Q>0$ and $\sigma\in (0,1]$.   
	Then, there exists  $C=C(n,\a,\s,Q)$ such that for all $E,F\subseteq \Rn$ with $\nu(E)<~\infty,\, |F|<~\infty$ and for all $f\in L^{\na}(F)$ with $\|f\|_\na\le 1$  we have
$$\int_E \exp\bigg[{\sigma\over |B_1|}|I_\a f(x)|^{n\over n-\alpha}\bigg]d\nu(x)\le C(\nu(E)+|F|^\sigma).\eqdef{18}$$
For given $E$ and $F$, if  there exists a ball $B(x_0,r_0)$ such that $|B(x_0,r_0)\cap F|=|B(x_0,r_0)|$, and $\nu(B(x_0,r)\cap E)\ge c_1 r^{\sigma n }$ for $r\le r_0$, with $c_1>0$,  then the 
exponential constant is sharp.
\par

The condition that $F$  contains a  ball (up to a set of zero measure) which has enough mass shared by $E$, is essentially necessary in order to guarantee sharpness in the above corollary. In general the sharp exponential  constant will depend on the relative geometry of the sets $E$ and $F$: the less  the mass they have in common, the larger  the sharp constant. This is a reflection of the fact that the potential becomes ``less effective" as the sets $E$ and $F$ get  more and more separated (in this regard, see [FM1, Remark 3, p. 5112]).

\pf Proof of Corollary 3.
In Corollary 2 take $N=\R^n$ endowed with the measure $\nu$ as in \eqref{b}, $M=\Rn$ with the Lebesgue measure,   $k(x,y)=|x-y|^{\a-n}$, $\b={n\over n-\alpha}$, so that under the assumption \eqref b we have (see also  [FM1, Lemma 9])
$$k_1^*(t)\le |B_1|^{{n-\a\over n}} t^{-{n-\a\over\ n}},\qquad k_2^*(t)\le Q^{n-\a\over\sigma n} t^{-{n-\a\over\sigma n}},\qquad t>0$$which implies \eqref{18}.

For the sharpness statement, assume that $F$ contains a ball of radius $r_0$ and center $x_0$, up to a set of zero measure. We can assume that $x_0=0$, and define for $0<\e<r_0$.
$$\phi_\e(y)=\cases{|y|^{-\a} & if $\e<|y|<r_0$\cr 0 & otherwise.\cr}$$
Then $\supp \phi_\e\subseteq F$, and 
$$\|\phi_\e\|_\na^\na=|B_1|\log{r_0^n\over \e^n}.$$
If $\tilde\phi_\e=\phi_\e/\|\phi_\e\|_\na$ then by routine  computations $|I_\a\tilde\phi_\e(x)|^\nna\ge |B_1|\log{r_0^n\over \e^n}-C$, provided $x\in B(0,\e/2)$. Assuming $\nu(E\cap B(0,\e/2)|\ge C \e^{\sigma n}$ we then get
$$\int_{E\cap B(0,\e/2)}\exp\bigg[\gamma|I_\a\tilde \phi_\e(x)|^{n\over n-\alpha}\bigg]d\nu(x)\ge C\e^{n(\sigma-\gamma|B_1|)}\to+\infty$$
if $\gamma>\sigma/|B_1|.$

\endpf

\smallskip
We now give the proof of Theorem 1. First let us note  that the estimate  given in \eqref{110}  will be  an immediate consequence of \eqref{109}, via  the following  elementary  lemma (see also [FM2, Lemma 9]):

\proclaim Lemma 4  (Exponential Regularization Lemma). Let  $(N,\nu)$ be a measure space and $1<p<\infty$, $\gamma>0$. Then for every $u\in L^p(N)$  we have
$$\int_{\{|u|\ge1\}} e^{\gamma|u|^{p'}}d\nu-e^\gamma\|u\|_p^p\le\int_{N}\exp_{[p-2]}\big[{\gamma |u|^{p'}}\big]d\nu\le\int_{\{|u|\ge1\}} e^{\gamma|u|^{p'}}d\nu+e^\gamma\|u\|_p^p.\eqdef{exp1}$$
In particular, the functional $\int_N \exp_{[p-2]}^{}\big[\gamma |u|^{p'}\big]$ is  bounded on a bounded subset $X$ of  $L^p$, if and only if $\int_{\{|u|\ge 1\}}\exp\big[\gamma |u|^{p'}\big]$ is  bounded  on $X$.
\smallskip
Moreover, for any $m=1,2,3....$ we have
$${\gamma^m\over m!}\|u\|_{mp'}^{mp'}\le \int_{N}\exp_{m-1}\big[{\gamma |u|^{p'}}\big]d\nu\le\int_{\{|u|\ge1\}} e^{\gamma|u|^{p'}}d\nu+e^\gamma\|u\chi_{|u|\le1}^{}\|_{mp'}^{mp'}.\eqdef{exp2}$$

\par

If the operator $T$ in Theorem 1 is also  continuous on $L^\bp(M)$, then the regularized exponential integral in \eqref{110} is clearly  bounded on the unit ball of $L^\bp(M)$.
For general ``well-behaved'' kernels $K$ continuity of $T$ in $L^\bp(M)$ implies critical integrability, in the form \eqref{108}  (see Proposition 2 below), but not viceversa. For example, consider on $\Rn$ $k(x,y)=K(x-y)$, where
$$K(x)=\cases{|x|^{\a-n}\qquad & if $|x|\le 1$\cr |x|^{\a-n-\delta} & if $|x|>1$\cr}$$
for any $\delta\in (0,\a)$. Here $\b=\nna,\;\bp=\nsa$. The kernel $K$  behaves like a Riesz kernel at $0$, and it satisfies the critical integrability condition. Yet, the convolution operator $K*f$ is not continuous on $L^\nsa(\Rn)$, since $K$ is nonnegative and not integrable.

\pf Proof of Theorem 1.
As we mentioned in the introduction, the proof of Theorem 1 is accomplished by making suitable  modifications and improvements  to the proof in [FM1, Theorem 1], in order to take into account the critical  integrability condition \eqref{108}, and by tracking down the various constants a little bit more carefully. For the convenience of the reader we will present here the beginning of the proof in enough details so that the role of \eqref{108} is highlighted, relegating the more technical parts  to the appendix.
\def\t{\tau}

Because of Lemma 4 it is enough to prove \eqref{109}. Indeed, 
$$\int_{N}\exp_{[\bp-2]}\bigg[{{\sigma\over A} |Tf(x)|^{\bp}}\bigg]d\nu(x)\le\int_{\{|Tf|\ge1\}} \exp\bigg[{\sigma\over A}|Tf(x)|^{\bp}\bigg]d\nu(x)+e^{\sigma\over A }\|Tf\|_\bp^\bp$$
and $|\{|Tf|\ge 1\}|\le \|Tf\|_\bp^\bp$. 

To prove \eqref{109} we show that for each $f\in L^\bp(M)$ the function $Tf$ is well-defined, finite a.e., and satisfies 
$$\int_0^{\nu(E)} \exp\bigg[{\s\over A}\big((Tf)^{**}(t)\big)^\beta\bigg]dt\le   Ce^{\s J_\t}\big(\nu(E)+\t^\s(1+J_\t)\big),\eqdef{0z}$$
where $ C=C(\b,\b_0,\gamma,H,A,B)$, under the hypotheses \eqref{106}, \eqref{107}, \eqref{108}, and with $$\int_0^\infty (f^*)^{\b'}\le 1.\eqdef{0d}$$
Below, and in the appendix, $C_j$ denotes a constant $\ge 1$, depending only on $A, B, \b,\s,p,H,\gamma$.
\smallskip
Without loss of generality we can assume that $k$ and $f$ are nonnegative. The first key element of the proof is the following improvement of Lemma 2 in [FM1] (which was itself an improvement of the original lemma due to O'Neil [ON, Lemma 1.5]):

\proclaim Lemma 5. Let $k: N\times M\to [0,\infty]$ be measurable and
$$k_1^*(t)\le D t^{-{1\over\beta}},\qquad k_2^*(t)\le B t^{-{1\over\sigma\beta}},\qquad t>0\eqdef{on1}$$
with $\beta>1$ and $0<\sigma\le 1$.
If  
$$\max\bigg\{1,{\b(1-\s),\over \b-1}\bigg\}\le p<{\b\over\b-1}=\b'\eqdef{on2}$$
and
$${1\over q}={1\over\s\b}+{1\over\s}\bigg({1\over p}-1\bigg),\qquad q>p\eqdef{on3}$$ then   there is a constant $C_0>0$ such that for any measurable $f:M\to [0,\infty]$ 
$$ (Tf)^{**}(t)\le C_0\,\max\big\{\tau^{-{\sigma\over  q }}, t^{-{1\over q}}\big\}\int_0^\tau u^{-1+{1\over p}}f^*(u) du+\esup\int_\tau^\infty k_1^*(x,u)f^*(u)du,\quad \forall t,\tau>0.\eqdef{on4}$$

\par \smallskip
\eject
\noin{\bf Notes.} \smallskip
\noin{\bf 1.} When $\sigma=1$ we can take $p=1$ and $q=\b$ in \eqref{on4}.\smallskip
\smallskip
\noin{\bf 2.} Even if $k$ is a kernel with arbitrary sign, conditions \eqref{on2} imply that $T$ is well defined on $L^p(M)$ and bounded from $L^p(M)$ to  $L^{q,\infty}(N)$ (see [A2] and our proof of Lemma 5 in the appendix).
\smallskip \noin{\bf 3.}
The earlier version of the lemma  given in  [FM1], and also used in [FM3],  has an inequality like \eqref{on4}, but with second term equal to $\int_\tau^\infty k_1^*(u)f^*(u)du$, which is larger than the one above.  With that version the conclusions of Theorem 1  can be proven under the stronger condition \eqref{stronger}.

\smallskip

 The proof of Lemma 5 is obtained by suitably modifying the proof given in [FM1]. The details are given in the appendix.
 \medskip
 
 Now we have  $k_1^*(t)\le A^{\bi} (1+H)t^{-\bi}$ for $t>0$, so that  by the improved O'Neil lemma above,  if $p$ is any fixed number satisfying \eqref{on2}, and $q$ is as \eqref{on3}, then
 for each $t>0$ 
$$(Tf)^{**}(t)\le C_0 t^{-{1\over q}}\int_0^{t^{1/\s}} u^{-1+{1\over p}}f^*(u) du+\esup \int_{t^{1/\s}}^\infty k_1^*(x,u)f^*(u)du.\eqdef{1}$$

\def\ts{\tau^\s}
If $t\ge \ts$
 then  \eqref{1} (combined into a single integral), H\"older's inequality and \eqref{0d} imply
$$\eqalign{(Tf)^{**}(t)& \le \bigg(\int_0^{t^{1/\s}} C_0^\b t^{-{\b\over q}}u^{\b\big({1\over p}-1\big) }du+\sup_{x\in N} \int_{t^{1/\s}}^\infty \big(k_1^*(x, u)\big)^\b du\bigg)^\bi\le\cr & \le \bigg( {C_0^\b\over\b\big({1\over p}-1\big)+1}+AJ_\t\bigg)^\bi.\cr}$$
and therefore
$$\int_{\ts}^{\nu(E)} \exp\bigg[{\s\over A}\big((Tf)^{**}(t)\big)^\beta\bigg]dt\le Ce^{{\s }J_\t} \big(\nu(E)- \ts\big)^+.$$\def\us{{1\over\s}}
On the interval $[0,\ts]$ unfortunately this simple argument fails  and we need to refine the more sophisticated analysis   in [A1] and [FM1].  From \eqref1 followed by the change of variables $t\to \ts t,\,u\to \tau u^\us$ we get
$$\eqalign{&\int_0^{ \ts} \exp\bigg[{\s\over A}\big(Tf)^{**}(t)\big)^\beta\bigg]dt\le \cr& \le \int_0^{ \ts} \exp\bigg[{\s\over A}\bigg(C_0 t^{-{1\over q}}\int_0^{t^{1/\s}} f^*(u)u^{-1+{1\over p}}du+\esup \int_{t^{1/\s}}^\infty k_1^*(x, u)f^*(u)du\bigg)^\b\,\bigg]dt\cr&\le \ts \int_0^1 \exp\bigg[{\s\over A}\bigg({C_0 t^{-{1\over q}}\t^{-{\s\over q}}\over\s}\int_0^{t} f^*(\t u^\us)\t^{1\over p}u^{-1+{1\over \s p}}du+\cr & \hskip14em +\us\esup  \int_{t}^\infty k_1^*(x, \t u^\us)f^*(\t u^\us)\t u^{-1+\us}du\bigg)^\b\,\bigg]dt\cr&
:=\ts \int_0^1 \exp\bigg[\bigg(C_2 t^{-{1\over q}}\int_0^t f_\t(u)u^{-{1\over\b}+{1\over q}}du+\esup \int_t^\infty k_\t(x, u)f_\t(u) du\bigg)^\b\,\bigg]dt
\cr}$$

where 
$$k_\t(x,u)=A^{-{1\over\b}}\t^{1\over \b}k_1^*(x,\t u^\us)u^{(-1+\us){1\over\b}},\qquad f_\t(u)=\s^{-{1\over\b'}}\t^{1\over\b'}f^*(\t u^\us)u^{(-1+\us){1\over\b'}}$$
and $C_2=C_0A^{-{1\over\beta}}$.
Note that 
$$\int_0^\infty f_\t(u)^{\b'}\le 1,\qquad \esup \int_\t^\infty k_\t(x,u)^\b du=\s J_\tau$$
and that 
$$\sup_{x\in N}k_\tau(x,u)\le u^{-{1\over\b}}\big(1+H(1+|\log(\t^\sigma u)|)^{-\gamma}\big),\qquad 0<u\le 1.$$

  Now make the further changes 
$$u=e^{-\xi},\quad t=e^{-\eta},\quad \phi(\xi)=f_\t(e^{- \xi}) e^{-{\xi\over\b'}}\eqdef{1q}$$
to obtain that 
$$\int_0^{\ts}\exp\bigg[{\s\over A}\big(Tf)^{**}(t)\big)^\beta\bigg]dt\le \ts \int_0^\infty e^{-F(\eta)}d\eta$$
where for each fixed $\eta\ge 0$
$$F(\eta)=\eta-\bigg(\sup_{x\in N}\int_{-\infty}^\infty g(x,\xi,\eta)\phi(\xi)d\xi\bigg)^\b$$
$$g(x,\xi,\eta)=\cases{k_\t(x,e^{-\xi})e^{-{\xi\over\b}} & if $\xi\le0$\cr
1+H(1+|\xi-\log\ts|)^{-\gamma} & if $0<\xi\le \eta$\cr
C_2e^{\eta-\xi\over q} & if $\eta<\xi$.\cr}$$
The next technical step is to run the Adams-Garsia machinery to prove that 
$$\int_0^\infty e^{-F(\eta)}d\eta\le C(1+J_\t)e^{\sigma J_\t}.\eqdef{AG}$$
The  details are given in the appendix.\endpf

\medskip

 For reasonable kernels the integrability condition \eqref{108} is essentially necessary in order to obtain exponential integrability or continuity from $L^\na$ to  any $L^p$:

\proclaim Theorem 6. Suppose that $\mu(M)=+\infty$ and that the kernel $k$ in \eqref{T}  satisfies 
$$k_1^*(t)<\infty,\qquad \forall t>0\eqdef{kstar00}$$
$$ \lim_{t\to+\infty} k_1^*(t)=0,\eqdef{kstar0}$$
and  that  for some sequence $x_m\in N$
$$\lim_{m\to+\infty}\int_1^\infty k_1^*(x_m,t)^\beta dt=+\infty.\eqdef{kstar1}$$
Suppose additionally  that for each $m\in \N$ there is a measurable set $B_m\subseteq N$  such that for some $\e_0>0$ and  all $\e>0$ small enough 
$$\mathop\int\limits_{\e<|k(x_m,y)|\le \e_0} \!\!\!k(x,y)k(x_m,y)|k(x_m,y)|^{\beta-2}d\mu(y)\ge C_0\!\!\!\!\!\!\mathop\int\limits_{\e<|k(x_m,y)|\le \e_0} \!\!\! |k(x_m,y)|^{\beta}d\mu(y),\qquad \forall x\in B_m\eqdef{reg}$$
for some $C_0>0$, independent of $\e$ (provided that the integral on the left hand side is well defined).
Then, there is a sequence of functions  $\{\Psi_m\}\subseteq L^\bp(M)$ with $\|\Psi_m\|_\bp=1$ such that 
\def\essinf{{\rm {ess\hskip.04in  inf}}}
$$\inf_{x\in B_m} |T \Psi_m(x)|\uparrow+\infty,\qquad \hbox{ as } \,m\to+\infty.$$
In particular if the $B_m$ can be chosen so that $\nu(B_m)\ge \nu_0>0$ for all $m$, then the operator $T$ defined in \eqref{T} cannot be continuous from  $L^\na$ to any $L^p$, and exponential integrability in the form \eqref{mt2} fails for any $\gamma>0$.
\par
\pf Proof. Conditions \eqref{kstar00}, \eqref{kstar0} guarantee that 
 $$\mu(\{y:|k(x,y)|=\infty\})=0,\qquad  \mu(\{y: |k(x,y)|>\e\})<\infty,\qquad \forall x\in N,\;\forall \e>0$$
 and that there is $T>0$ such that  
 $$\mu(\{y:\;|k(x,y)|>\e_0\}\le T,\qquad \forall x\in N.$$
Also, \eqref{kstar00} implies that \eqref{kstar1} holds if and only if for any $\delta>0$
$$\lim_{m\to+\infty}\int_\delta^\infty k_1^*(x_m,t)^\beta dt=+\infty.$$

The level sets 
$$F_\e^m=\{y\in M:\;|k(x_m,y)|>\e\},\qquad \e\ge0$$ satisfy, for each given $m$,  
$0<\mu(F_\e^m)\uparrow\mu(F_0^m)$ as $\e\downarrow 0$,  and $\mu(F_{\e_0}^m)\le T$, some $T>0$. Hence, for fixed $m$ and as $\e\downarrow 0$
$$\eqalign{\Gamma_\e^m:&=\int_{F_\e^m\setminus F_{\e_0}^m}|k(x_m,y)|^\beta d\mu(y)=\int_{\mu(F_{\e_0}^m)}^{\mu(F_\e^m)}k_1^*(x_m,t)^\beta dt\uparrow\int_{\mu(F_{\e_0}^m)}^{\mu(F_0^m)}k_1^*(x_m,t)^\beta dt=\cr& =\int_{\mu(F_{\e_0}^m)}^\infty k^*(x_m,t)^\beta dt\ge  \int_{T}^\infty k_1^*(x_m,t)^\beta dt\cr} $$
Using  \eqref{kstar1} we can find subsequences $m_k\uparrow\infty$ and $\e_k\downarrow 0$ such that $\Gamma_{\e_k}^{m_k}\to+\infty$ as $k\to+\infty$. By replacing the original sequence $\{x_m\}$ with $\{x_{m_k}\}$ we may  assume that $\Gamma_{\e_m}^m\to+\infty$ for a sequence $x_m$ and $\e_m\downarrow 0$.

If we let
$$\Phi_m(y)=k(x_m,y)|k(x_m,y)|^{\beta-2}\chi_{ F_{\e_m^{}}^m\setminus F_{\e_0}^m}^{}$$
then
$$\|\Phi_m\|_\bp^\bp= \Gamma_{\e_m}^m$$
and, using the hypothesis, $T\Phi_m(x)$ is well-defined on some $B_m$ with
$$T\Phi_m(x)=\int_{F_{\e_m}^m\setminus F_{\e_0}^m }k(x,y)k(x_m,y)|k(x_m,y)|^{\beta-2} d\mu(y)\ge C_0 \Gamma_{\e_m}^m,\qquad \forall x\in B_m.$$

The result follows upon taking $\Psi_m=\Phi_m\|\Phi_m\|_\bp^{-1}$.


\endpf

\noin{\bf Note.} If one of the conditions \eqref{kstar00}, \eqref{kstar0} is not satisfied (regardless of \eqref{kstar1}), then the conclusion of the theorem still holds, provided that $\mu$ is semifinite and \eqref{reg} holds.
\medskip

Certainly condition \eqref{kstar1} is met for some sequence $\{x_m\}$ if \eqref{108} fails. It does seem unavoidable to impose some sort of ``regularity''  condition equivalent to \eqref{108}, however. On $\Rn$, for example, condition \eqref{reg} is typically verified when $B_m$ is a small ball around $x_m$ and for $x\in B_m$, $k(x,y)-k(x_m,y)$ decays better than $k(x_m,y)$ as $y\to+\infty$. For  convolution kernels which are radially decreasing (even just outside a large ball) $k(x,y)=K(|x-y|)\ge0$  condition \eqref{reg}  is easily verified, since both $K(|x-y|)$ and $K(|y|)$ are greater than $K(2|y|)$, for $|x|\le 1$ and $|y|\ge2$ 
\bigskip
\centerline{\scaps 3. Riesz subcritical kernels and domains  in $\Rn$}\bigskip

The most obvious application of Theorems 1 and 6 is to convolution operators on $\Rn$:

\proclaim Theorem 7. Let $0<\a<n$,  and suppose that $K:\Rn\setminus\{0\}\to\R$ is measurable and satisfies the conditions
$$K(x)=g(x^*)|x|^{\a-n}+ O(|x|^{\a-n+\e}),\qquad |x|\le R,\; x^*={x\over|x|}\eqdef{111}$$
$$K\in L^{\nna}\cap L^\infty\big( \{|x|\ge R\}\big)\eqdef{112}$$
for some $\e,R>0$, where $g\in L^{\nna}(S^{n-1})$. Then $K*f$ is finite a.e for $f\in L^{\nsa}(\Rn)$, and   there exists $C>0$ such that for all $f\in L^{\nsa}(\Rn)$ with $\|f\|_\na\le 1$, and for each measurable $E\subseteq \Rn$ with $|E|<\infty$
$$\int_E\exp\bigg[{1\over A_g}|K*f(x)|^\nna\bigg]dx\le  C(1+|E|),\eqdef{112a}$$
where $$A_g={{1\over n}{\isn|g(\omega)|^{\nna}d\omega}}.\eqdef{114}$$
Moreover, 
$$\irn\exp_{[\nsa-2]}\bigg[{1\over A_g}|K*f(x)|^\nna\bigg]dx\le  C\big(1+\|K *f\|_\na^\na\big).\eqdef{113}$$
for all $f\in L^\nsa(\Rn)$ such that $\|f\|_\na\le 1$.
If $g$ is smooth, then the exponential constant in \eqref{112a} (if $|E|>0$) and in \eqref{113} is sharp.\par

Note that the ``big O" notation  in \eqref{111} means that  $|O(|x|^{\a-n+\e})|\le C|x|^{\a-n+\e}$, for $|x|\le R$. 

\smallskip
\pf Proof of Theorem 7. In Theorem 1 let $N=E$ and $M=\Rn$ with the Lebesgue measure, and let $k(x,y)=K(x-y)$, $\beta=\nna$. Note first that $k_1^*(x,t)=k_1^*(t)=k_2^*(t)$, for all $x,t$.  Let $\tau=|B(0,R)|$ and suppose that $|K(x)|\le M$ when $|x|\ge R$. If $\lambda(s)=|\{x\in \R^n: |K(x)|>s\}|$ denotes the distribution function of $K$ relative to $\R^n$, then for $s>M$ the distribution function of $K$ relative to the ball of radius $R$ coincides with $\lambda(s)$. This means that $k_1^*(t)$ is the same as the corresponding rearrangement relative to  $B(0,R)$, when $t<\tau$. The proof that \eqref{111} implies \eqref{106} and  \eqref{107} for small $t$, on sets of finite measure and therefore for $t\in (0,\tau]$,  has been done in [FM1, Lemma 9]. Note that the proof there was done in the case $g$ bounded on the sphere, but it works even in our more general hypothesis.

It is enough to check that  \eqref{112} implies   \eqref{108}    (from which  \eqref{107} follows for all $t$, since $k_2^*=k_1^*$ is decreasing and finite).
The proof of this fact is  straightforward. Let $|K(x)|\le M$ for $|x|\ge R$, and let 
$$\wtilde K(x)=\cases{|K(x)| & if $|x|\ge R$\cr M & if $|x|<R.$\cr}$$
If $\lambda(s)$ and $\wtilde\lambda(s)$ denote the distribution functions of $K,\wtilde K$ respectively, then $\wtilde\lambda(s)\ge \lambda(s)$ for $s< M$, and $\wtilde\lambda(M^-)\ge |B_1|R^n$. Hence, if $ k_1^*,\,\wtilde k_1^*$ denote the rearrangements of $ K,\, \wtilde K$ resp., then $\wtilde k_1^*(t)\ge k_1^*(t)$ for $t\ge |B_1|R^n$.
Obviously,  $\wtilde K\in L^\nna(\Rn)$, so \eqref{108} follows with $\tau =|B_1|R^n$.
This proves  inequality \eqref{112a}, and therefore \eqref{113}.

The proof of the sharpness statement is the same as that of [FM1, Theorem 8]. In particular, assuming, without loss of generality, that $|E\cap B(0,\e)|\ge c_0 \e^n$ for small $\e$ and for some $c_0>0$,  one can take the extremal family of functions 
$$\phi_\e(y)=\cases{K(y)|K(y)|^{\nna-2} , & if $\e<|y|<1$\cr 0 & otherwise\cr}\eqdef{extremals}$$
and show that along the normalized family $\phi_\e/\|\phi_\e\|_\na$ the exponential integrals in \eqref{112a} and \eqref{113} are saturated.
 \endpf

 
Generally speaking,  non-homogeneous invertible elliptic  operators will have  kernels satisfying \eqref{112}, and for those operators a sharp Moser-Trudinger inequality will hold.  As a first example, consider the Bessel potential $(I-\Delta)^{\a\over2}$, whose fundamental solution $G_\a(x)$ behaves like the Riesz potential locally, and decays exponentially at infinity. In fact, the aforementioned results by [A, Thm. 3],  [LL2] and [RS] are immediate consequences of Theorem 7, and the fact that $\|G_\a* f\|_p\le \|f\|_p$:

\proclaim Theorem 8. If $0<\a<n$ then there exists $C$ such that for all $u\in W^{\a,\nsa}(\Rn)$ so that 
$$\|(I-\Delta)^{\a\over2}u\|_\nsa\le 1$$
we have
$$\irn\exp_{[\nsa-2]}\bigg[\gamma\big((-\Delta)^{\am}\big)|u(x)|^\nna\bigg]dx\le  C.\eqdef{bessel}$$
and the exponential constant is sharp.\par

In this paper we  define $W^{\a,\nsa}(\Rn)$ for $\a>0,p>1$ to be the space of Bessel potentials:
$$\eqalign{W^{\a,p}(\Rn)&=\{u\in L^p(\Rn):\;(I-\Delta)^\am u\in L^p(\Rn)\}=\{G_\a*f,\,f\in L^p(\Rn)\}\cr&=\{u\in L^p(\R^n):\;(-\Delta)^\am u\in L^p(\Rn)\}.\cr}\eqdef{W}$$

In section 4, Theorem 11, we will exhibit  a class of  non-homogeneous, elliptic, invertible, linear partial differential  operators $P$  whose inverses have   kernels satisfying \eqref{111} and \eqref{112}, and therefore a sharp Moser-Trudinger inequality holds for such $P$ (Theorem 14).

\medskip
We point out that  Theorem 7 can be   formulated so as to accommodate  more general (non-convolution) kernels satisfying
$$K(x,y)=g(x,(x-y)^*)|x-y|^{\a-n}+O(|x-y|^{\a-n+\e}),$$
together with suitable integrability and boundedness conditions at infinity, in the same spirit as in [FM1, Thm. 8].

We also remark  that Theorem 7 could have been stated in the slightly more general situation  where the convolution  operator is acting on $L^\nsa(\Omega)$, where $\Omega\subseteq \Rn$ is an arbitrary measurable set of $\Rn$. In this case the conclusion holds provided that $K\in L^\infty\cap L^\nna(\{x\in \Omega:|x|\ge R\})$. We find that the latter condition is of little applicability if $\Omega\neq \Rn$, in which case one is better off checking out the corresponding critical integrability condition on $k_1^*(x,t)$, the rearrangement of $K(x-\cdot)$ with respect to $\Omega$. 

When $K$ is the Riesz kernel, however, \eqref{108} leads to an interesting geometric  condition on $\Omega$, under which inequality \eqref{112a} holds under the sole condition that $\|f\|_{L^\na(\Omega)}^{}\le 1$, as expressed in the next theorem.
\def\G{\Lambda_\Omega^{}}
\medskip For a measurable  set $\Omega\subseteq \Rn$ define for $r>0$ 
$$\G(x,r)={|\Omega\cap B(x,r)|},\qquad x\in \Omega,$$
$$\G(r)=\sup_{x\in\Omega} \G(x,r).\eqdef{supG}$$

\proclaim Theorem 9. Let $\Omega\subseteq\Rn$ be measurable and such that 
$$\sup_{x\in \Omega}\,\int_1^\infty {\G(x,r)\over r^{n+1}}dr<\infty.\eqdef G$$
   Then, for  $0<\a<n$, there exists $C>0$ such that for all $f\in L^{\nsa}(\Omega)$ with $\|f\|_\na\le 1$, and for each measurable $E\subseteq \Omega$ with $|E|<\infty$
$$\int_E\exp\bigg[{1\over |B_1|}|I_\a f(x)|^\nna\bigg]dx\le  C(1+|E|).\eqdef{G1}$$
Moreover, 
$$\int_\Omega\exp_{[\nsa-2]}\bigg[{1\over |B_1|}|I_\a f (x)|^\nna\bigg]dx\le  C\big(1+\|I_\a f\|_{L^\na(\Omega)}^\na\big).\eqdef{G2}$$
for all $f\in L^\nsa(\Omega)$ such that $\|f\|_\na\le 1$.
If there is $x_0\in\Omega$ and $r_0>0$ such that $|\Omega\cap B(x_0,r_0)|=|B(x_0,r_0)|$, and if $|E|>0$,  then the exponential constant in \eqref{G1} and \eqref{G2} is sharp.  \smallskip 
Conversely, if \eqref G is not satisfied then \eqref{G1} cannot hold, in fact there is a sequence of functions $\Psi_m\in L^\na(\Omega)$, with $\|\Psi_m\|_\na=1$ and $r_1,\delta>0$ such that with $B_m=B(x_m,r_1)\cap \Omega$  and $|B_m|\ge\delta$ we have 
$$\inf_{x\in B_m}|I_\a\Psi_m(x)|\to+\infty.$$\par

Condition  \eqref{G} is independent of $\alpha$, and expresses the Riesz subcriticality of the Riesz potential restricted to the measurable set  $\Omega$.

We also note  that \eqref{G} is implied by the stronger condition
$$\int_1^\infty{\G(r)\over r^{n+1}}<\infty\eqdef{stronger}$$
with $\G(r)$ defined in \eqref{supG}. As the following example shows, it is possible to  construct an $\Omega$ such that \eqref{G} (and hence \eqref{G1}and  \eqref{G2}) holds but \eqref{stronger} fails. 
\def\k{{\bf k}}\def\ee{{\bf e}} Let $c_m,\,R_m,\, \e_m$ be positive sequences with $\e_m<\half$, $c_m\uparrow+\infty,\;R_m\uparrow+\infty,$  $c_m-c_{m-1}>>R_m$.
If $k=(k_1,...,k_n)\in \Z_+^n$ and $x_m=c_m\ee_1$,  let 
$$\Omega_m=\bigcup_{ |k-x_m|<R_m} B(k,\e_m),\qquad \Omega_*=\bigcup_{m=1}^\infty \Omega_m.\eqdef{dominio1}$$
In other words, $\Omega_m$ is the union of all the balls of radius $\e_m$ centered at the integer points contained in the ball of radius $R_m$ and center $x_m$.
\def\Gstar{{\Lambda_{\Omega_*}}}
With this in mind, it is possible to  show  that there is $C$ independent of $m$ such that for all $x\in \Omega_m$
$$\int_{1}^\infty{\Gstar(x,r)\over r^{n+1}}dr\le C\bigg[\e_m^n\log R_m+\sum_{\ell=m}^\infty\ell\bigg( {\e_{\ell}R_{\ell}\over c_{\ell}}\bigg)^{\!\! n}\,\bigg]\eqdef{supint1}$$
and 
$$\int_{1}^\infty{\Gstar(r)\over r^{n+1}}dr\ge C\sum_{m=2}^\infty \e_m^n\log{R_{m}\over R_{m-1}}\eqdef{supint2}$$

Choosing  for example
$$R_m=4^m,\qquad \e_m^n\log R_m=2^{-n},\qquad c_m=100^m$$
we have that the series in \eqref{supint1} is finite and the one in   \eqref{supint2} is infinite. Note that the Adams inequality for $L^\nsa(\Omega_*)$ is guaranteed by Theorem 9, whereas it would be hard to determine this fact using previously known methods.

\smallskip\smallskip
 \pf Proof of Theorem 9. Let us apply Theorem 1 with $M=N=\Omega$ and $\mu=\nu=$ Lebesgue measure. Since $I_\a$ has kernel $|x-y|^{\a-n}$ we have $k_1^*=k_2^*$, where rearrangement is with respect to $\Omega$. For simplicity let us drop the index ``1": and let  
$$\lambda(x,s)=|\{y\in\Omega:\,|x-y|^{\a-n}>s\}|=\big|\Omega\cap B\big(x,s^{-{1\over n-\a}}\big)\big|=\G\big(x,s^{-{1\over n-\a}}\big)$$
  $$k^*(x,t)=\inf\{s>0:\,\lambda(x,s)\le t\},\qquad$$
$$\lambda(s)=\sup_{x\in\Omega} \lambda(x,s),\qquad k^*(t)=\sup_{x\in\Omega} k^*(x,t).$$

Observe here that

$$\lambda(s)=|\{t:k_1^*(t)>s\}|\eqdef{lambda}$$
from which it follows
$$k_1^*(t)=\inf\{s>0:\,\lambda(s)\le t\}.\eqdef{lambda1}$$
(Equation \eqref{lambda} holds for general kernels $k$, and it was stated in [FM1, p. 5073, ``Fact~3"]. The proof is based on the fact that for each  $s>0$ and each $x\in \Omega$ there exists $T_{x,s}>0$ such that 
$\{t:k_1^*(x,t)>s\}=(0,T_{x,s})$.)

From general facts about rearrangements we have that $\lambda(x,\cdot),k^*(x,\cdot)$ are decreasing and right-continuous. Moreover 
$k^*(x,\lambda(x,s))\le s$ and $\lambda(x, k^*(x,t))\le t$. However, in this case we have also that $\lambda(x,s)$ is actually continuous  (in fact locally Lipschitz) in $s$, hence $\lambda(x, k^*(x,t))=t$ for all $t>0$, and all $x\in\Omega$.

Obviously,
$\G(x,r)\le |B_1|r^n $ for every $x\in \Omega$, and $r>0$,  hence 
$$k^*(t)\le |B_1|^{{n-\a\over n}}t^{-{n-\a\over n}},\qquad \forall t>0\eqdef{stimak1}$$ and condition \eqref{106} is verified with $A=|B_1|$ and $\b=\nna$.

Now note that if $\phi:(0,\infty)\to[0,\infty)$ is decreasing and right-continuous, and if $\lambda_\phi(s)$ denote its distribution function, then for each $t_0>0$
$$\int_0^{\phi(t_0)}\lambda_\phi(s) s^{p-1}ds={1\over   p}\int_{t_0}^\infty \phi(t)^p dt+t_0 \phi(t_0)^{p-1}\eqdef{formula}$$
(this is a consequence of Fubini's theorem).
Hence, applying this to $\phi(t)=k^*(x,t)$, $p=\nna$,  $t_0=\lambda(x,1)$ we get 
$$\eqalign{\int_{\lambda(1)}^\infty\big( k^*(x,t)\big)^\nna dt\le \int_{\lambda(x,1)}^\infty\big( k^*(x,t)\big)^\nna dt&\le{n\over n-\a}\int_0^{k^*(x,\lambda(x,1))}\lambda(s) s^{\nna-1}ds\cr&\hskip-3em  \le{n\over n-\a}\int_0^1\lambda(x,s) s^{\a \over n-\a}ds  = n \int_1^\infty {\G(x,r)\over r^{n+1}} dr\cr}\eqdef{G10}$$
which implies that condition \eqref{108} holds for any $\tau>\lambda(1)$.

This shows that the conditions of Theorem 1 are met, and the exponential inequalities follow. The sharpness of the exponential constant follows exactly as in the proof of Lemma~3.

Conversely, suppose that \eqref G does not hold. Then we can find a sequence $\{x_m\}\subseteq \Omega$ such that 

$$\int_1^\infty {\G(x_m,r)\over r^{n+1}}dr\to+\infty.\eqdef {G8}$$

We would like to apply Theorem 6. Conditions \eqref{kstar00}, \eqref{kstar0} of Theorem 6 are a consequence of \eqref{stimak1}.  To show \eqref{kstar1},
 apply formula \eqref{formula}, with $\phi(t)=k^*(x_m,t)$ and $t_0=1$, together with \eqref{stimak1} 
$$\eqalign{\int_1^\infty\big( k^*(x_m,t)\big)^\nna dt&={n\over n-\a}\int_0^{\d_m}\lambda(s) s^{\a\over n-\a}ds-\nna\d_m^{\a\over n-\a}\cr& \ge n\int_{\d_m^{-{1\over n-\a}}}^\infty {\G(x_m,r)\over r^{n+1}} dr-\nna|B_1|^{\a/n}\cr}$$
where we let
$$\d_m=k^*(x_m,1).$$
If $\d_m\ge1$ for infinitely many $m$, then for a subsequence $\{x_{m_k}\}$
$$\int_1^\infty\big( k^*(x_{m_k},t)\big)^\nna dt\ge n\int_1^\infty {\G(x_{m_k},r)\over r^{n+1}} dr-n|B_1|^{\a/n}\to+\infty$$
so condition \eqref{kstar1} is satisfied up to passing to a subsequence.

If instead $\d_m=k^*(x_m,1)<1$ for all $m$ large enough, since $\lambda(x_m,\d_m)=1$ (owing to the continuity of $\lambda(x,\cdot)$) we have $\G(x_m,\d_m^{-{1\over n-\a}})=1$ and therefore (since $\G(x_m,\cdot)$ is increasing) 
$$\int_1^{\d_m^{-{1\over n-\a}} }{\G(x_m,r)\over r^{n+1}}dr\le \int_1^\infty {1\over r^{n+1}}={1\over n}$$
$$\int_1^\infty\big( k^*(x_{m_k},t)\big)^\nna dt\ge n\int_1^\infty {\G(x_{m_k},r)\over r^{n+1}} dr-1-n|B_1|^{\a/n}\to+\infty$$
and condition \eqref{kstar1} is satisfied even in this case.
\def\tk{\tilde k}

 To verify  the integral condition in \eqref{reg} we proceed as follows.  Define 
 $$\tilde k(t)=\sup_m k^*(x_m,t),\qquad \tilde \lambda (s)=\sup_m \lambda(x_m,s).$$ Clearly the relations in \eqref{lambda}, \eqref{lambda1} continue to hold with $\tilde\lambda,\, \tk$ in place of $\lambda, \,k_1^*$, and  \eqref{kstar1}  implies $\tk(t)>0$ for all $t$. Pick any $t_1>0$. Since $\tk(t)\to 0$ as $t\to +\infty$ (from \eqref{stimak1}), we can find $t_0>t_1$ such that $0<\tk(t_0)<\tk(t_1)$. Let
 $$\e_0=\tk(t_0),\quad \e_1=\tk(t_1),\quad r_0=\e_0^{-{1\over n-\a}}>r_1=\e_1^{-{1\over n-\a}}$$
 With this choice of $\e_0$  the inequality in  \eqref{kstar1} can be written as
   $$\mathop\int\limits_{\{y\in\Omega:\,r_0\le |y-x_m|<\e^{-{1\over n-\a}}\}}\hskip-2em |x-y|^{\a-n}|y-x_m|^{-\a}dy\ge C_0\hskip-3em\mathop\int\limits_{\{y\in\Omega:\,r_0\le |y-x_m|<\e^{-{1\over n-\a}}\}}\hskip-2em |y-x_m|^{-n}dy.$$
The inequality is certainly verified if $x$ belongs to the set
$$B_m=\{x\in\Omega:\, |x-x_m|<r_1\}$$
since $|y-x|\le |y-x_m|+r_1\le 2|y-x_m|$, if $|y-x_m|\ge r_0,\;|x-x_m|<r_1$.

Moreover,  $|B_m|=\lambda(x_m,\e_1)$, and since $\tilde \lambda(\e_1)=t_1$, we can pass to a subsequence of $x_m$ in order to guarantee that $|B_m|\ge t_1/2$, for all $m$.
\endpf
\noin{\bf Remark.}  Under the hypothesis \eqref{G}, estimate \eqref{G1} actually holds for {\it any} measurable set $E\subseteq \R^n$. To see this, note that for each fixed $x\in\R^n$ if $\Omega$ has positive measure then $\G(x,r)$ is increasing in $r$ and eventually positive. If  $R_x:=\inf\{r>1: \G(x,r)>~0\}$ then  $\G(x,r)\le \G(y,3r)$
for any $y\in A_x:=B(x,2R_x)\cap \Omega$ (of positive measure) and for any $r\ge R_x$. Hence for all $y\in A_x$
$$\int_1^\infty {\G(x,r)\over r^{n+1}}dr\le\int_{R_x}^\infty {\G(y,3r)\over r^{n+1}}dr\le 3^n\sup_{y\in \Omega} \int_1^\infty {\G(y,r)\over r^{n+1}}dr$$
which implies that $$\sup_{x\in\sR^n}\int_1^\infty {\G(x,r)\over r^{n+1}}dr<\infty,$$
and one can  proceed as in the proof of Theorem 9 with $N=\Rn$  instead of $N=\Omega$.

Theorem 9  gives a sufficient condition independent of $\alpha$, i.e. \eqref{G}, under which  Adams' original result \eqref{mt2}  holds for domains of infinite measure. In a sense, the condition says that the domain $\Omega$ misses enough dimensions at infinity. Examples of such domains are ``strips" namely
$\Omega=\{x\in\R^n: x_j\in[a_j,b_j], j=1,...,k\}$, ($k<n$), in which case it is easy to see that $\G(r)\sim C r^{n-k}$  as $r\to+\infty$.   
It is not hard to construct domains of infinite measure so that the corresponding  $\G(r)$ has  prescribed order of growth, within the upper bound $|B_1|r^n$. Take any smooth $h:[0,\infty)\to[0,\infty)$, $h(0)=0$,  strictly increasing to $+\infty$ and with $h(x+1)-h(x)$   increasing (for example $h$   convex), and   let\def\m{{\bf m}} 
$$\Omega=\bigcup_{\m\in\sZ_+^n}^\infty B(C_\m,\delta_0),\qquad\m=(m_1,...,m_n)\in \Z_+^n,\;\,\quad  C_\m=\big(h(m_1),...,h(m_n)\big)\eqdef{dominio}$$
where $\delta_0>0$ is chosen so that all $B(C_\m,10\d_0)$ are pairwise disjoint. Then, one can check that for some $c_1,c_2,c_3>0$
$$ c_1\Big(h^{-1}\Big({r \over\sqrt n}\Big)\Big)^n \le \G(0,r)\le c_2  \big(h^{-1}(r+1)\big)^n,\qquad  r>r_0,\eqdef{stimaG1}$$
and 
$$\G(x,r)\le c_3 \G(0,r\sqrt n),\qquad x\in \Omega,\;r\ge r_0\eqdef{stimaG2}$$
some $r_0$ large enough. 
For the details of the proof see the Appendix.

Estimates \eqref{stimaG1} and \eqref{stimaG2} give
$$ c_1\Big(h^{-1}\Big({r \over\sqrt n}\Big)\Big)^n \le \G(r)\le c_4  \big(h^{-1}(r\sqrt n+1)\big)^n$$
for all $r$ large enough. With this in mind, one can, for example,  produce an $h$ as above so that  $\G(r)$ grows like $r^{n-\delta} (\log r)^{-q}$, for large $r$, for any $\delta,q$ with $0\le \delta<n,\, q\ge0$.
\smallskip

\medskip
\proclaim Corollary 10. If $\Omega$ is open and Riesz subcritical, i.e. if  it satisfies \eqref{G}, then for each $\alpha$ integer in $(0,n)$ there is $C$ such that  for all $u\in W_0^{\a,\na}(\Omega)$ with $\|\nabla^\a u\|_\na\le1$ we have
$$\int_{\Omega}\exp_{[\nsa-2]}\Big[\gamma(\nabla^\a)|u(x)|^\nna\Big]dx\le C(1+\|u\|_{L^\na(\Omega)}^\na)\eqdef{MTS0}$$
and the exponential constant is sharp.\par
 
 From the above estimate it is clear that if the Poincar\'e inequality $\|u\|_\na\le C\|\nabla^\a u\|_\na$ holds in $W_0^{\a,\na}(\Omega)$, then there is uniformity on the right hand side of \eqref{MTS0}. We shall return to this connection with the Poincar\'e inequality in Section 6.
 \bigskip 

\bigskip \centerline{\scaps 4. Riesz subcritical  fundamental solutions of elliptic differential operators  }
\centerline{\scaps with constant coefficients on $\Rn$}

\bigskip

Let  us consider an elliptic differential operator of  order $\alpha<n$ with constant complex  coefficients and acting on $C^\infty_c(\Rn)$
$$Pu=\sum_{|k|\le\a}a_k D^k u$$
where $k=(k_1,...,k_n)$ denotes a nonnegative multiindex in $\Z^n$. We will let
$$p(\xi):=P(2\pi i \xi)=\sum_{|k|\le\a}a_k (2\pi i\xi)^k$$
and define the strictly homogeneous principal symbol of $P$ as
$$p_\a(\xi):=P_\a(2\pi i\xi):=(2\pi i)^\a\sum_{|k|=\a}a_k \xi^k.$$

For simplicity  we only consider the case $a_k\in \R$, in which case $\a$ is even and ``$P$ elliptic"  means that
$$|p_\a(\xi)|\ge c_0|\xi|^\a,\qquad \xi\in \R^n\eqdef{ellip0}$$
for some $c_0>0$.

\def\k{{\kappa}}
\def\kp{K_P^{}}
 It is well known that $P$ has a fundamental solution,  given by a function ${K_P^{}}^{}$  which is $C^\infty$ outside the origin, and which is  formally the inverse Fourier transform of $1/p( \xi)$ i.e.
 $${K_P^{}}^{}(x)=\irn {e^{2\pi i x\cdot\xi}\over p(\xi)}\,d\xi.$$
 With this notation we have that
 $$u=\kp*(Pu)\qquad P(\kp*f)=f.$$
 
In what follow we will consider the case $p(x)\neq0$ for $x\neq0$, so that   $\kp$ has a singularity only at $0$. Indeed, 
 a  formula using classical integrals for ${K_P^{}}^{}$ can be written for example as follows: 
$${K_P^{}}^{}(x)=\irn {\eta(\xi)\over p(\xi)}\,e^{2\pi i x\cdot\xi}d\xi+{1 \over(2\pi|x|)^{2\ell}}\irn \Delta^\ell\bigg({1-\eta(\xi)\over p(\xi)}\bigg)\,e^ {2\pi i x\cdot\xi}d\xi,\eqdef {K1}$$
for $x\neq0$, were $\eta$ is a smooth cutoff which is 1 for $|x|\le 1$ and 0 for $|x|\ge 2$, for any $\ell>{n-\a\over2}$. This follows by writing $e^{2\pi ix\cdot\xi}=(2\pi |x|)^{-2\ell}\Delta^\ell\big(e^{2\pi ix\cdot\xi}\big)$ and integrating by parts.

The first term in \eqref{K1} is in $C^\infty(\Rn)$ and controls the decay of $K_P$ at infinity. The second term is $C^\infty(\Rn\setminus 0)$ and controls the singularity  of $K_P$ at the origin. To analyze the singularity,    write $p=p_\a+p_m$, for some $p_m$ polynomial of order $m<\a$ and 
$${1\over p_\a+p_m}=\sum_{k=0}^\ell (-1)^k{p_m^k\over p_\a^{k+1}}+{(-1)^{\ell+1}p_m^{\ell+1}\over p_\a^{\ell+1}(p_\a+p_m)}.\eqdef{expan}$$
The last term above is integrable outside a ball if $\ell$ large, whereas the other terms can be arranged  into a finite sum, where the first term is $1/p_\a$ and the other terms are all homogeneous of order $<-\a$. From this one obtains that
$$\kp(x)={g_P^{}}^{}(x^*)|x|^{\a-n}+O(|x|^{\a-n+\e}),\qquad |x|\le1\eqdef{K2}$$
where $x^*=x/|x|$, and where 
$${g_P^{}}^{}(x)=\irn {e^{2\pi i x\cdot\xi}\over p_\a(\xi)}\,d\xi\eqdef{gp}$$
in the sense of distributions
 (see also [FM1], formulas (67), (69)).  This is precisely the local asymptotic expansion \eqref{111} of Theorem 7, which  has already been used in a more general context in [FM1] to prove the sharp Moser-Trudinger inequality for $P$ on bounded domains. 
The validity of the sharp Moser-Trudinger inequality for $P$ in the form \eqref{mt1} on the whole of $\Rn$ is therefore relying on the critical integrability of $K_P$ at infinity. 

\smallskip\noin{\bf{Question:}} {\it Which non-homogeneous elliptic differential operators with constant coefficients on $\Rn$ have a fundamental solution  which is Riesz subcritical?}

\smallskip
Note that if $P$ is homogeneous of order $\alpha$ then $K_P$ is homogeneous of order $\a-n$ and critical integrability does not hold.
\smallskip

The precise asymptotic behavior of $\kp$ for large value of $x$ is not so obvious to figure out. It is well known that if $p$ never vanishes then $\kp$ decays exponentially at infinity. For the case  $p(0)=0$ not much seems to be known in the literature other than a few special cases. From \eqref{K1} we  see that  $\kp(x)\to0$ as $x\to\infty$, and in particular  it is bounded outside a ball centered at 0.

If the lowest order terms of $P$ form an elliptic homogeneous operator of order $m<\alpha$,  then one can show that  $|\kp(x)|\le C |x|^{m-n}$ for large $x$, so critical integrability  holds.  This can be seen by  writing $p=p_\a'+p_m'$, for some $p_m'$ elliptic homogeneous or order $m$,  $p_\a'$ is elliptic of order $\a$, and   using  formula \eqref{expan} with $p_m'$ in place of $p_\a$ and $p_\a'$ in place of $p_m$.

If $P=\big(\nabla^T\A\nabla+b^T\nabla\big)^{\a/2}$, where  $\A$ is a real, symmetric and positive definite matrix and $b\in \R^n\setminus\zero$, then the fundamental solution $\kp$ can be explicitly computed via linear transformations  from the one for the Bessel operator (see [OW, (2.5.3)] and [Lo], formula (7)). For example, when $\A=I$, the identity matrix,  we have
$$\kp(x)={(-1)^{\am}\over 2^{n-1}\pi^{n/2}\Gamma\big(\am\big)}{|x|^{\a-n\over2}\over|b|^{\a-n\over2}} e^{-{b\cdot x\over2}} K_{n-\a\over2}\big(\half|x||b|\big)\eqdef{FS1}$$
where $K_\nu$ is the modified Bessel function of the second kind. From this formula one obtains that for  $|x|\ge R/|b|$  we have $|\kp(x)|\le C_b |x|^{\a-n-1\over2}e^{-{|x||b|+b\cdot x\over2}}$, and a little computation in polar coordinates reveals that $\kp$ is in $L^p(|x|\ge1)$ for 
$p>{n+1\over n+1-\a}$, which  includes the case $p={n\over n-\a}$. For general $\A$ critical integrability follows for the same reasons, since the fundamental solution is obtained from \eqref{FS1} by multiplying it by $(\det\A)^{-\half}$ and by replacing $x$ with $\A^{-1/2}x$ and $b$ with $\A^{-1/2}b$.
\smallskip
For  $\a\ge n/2$ we can prove the following result:
\smallskip
\eject 
\proclaim Theorem 11. If $P$ is a non-homogeneous  elliptic differential operator  with constant coefficients of order $\a$, with ${n\over2}\le \a<n$ and such that for some $c_1>0$ 
$$|p(\xi)|\ge c_1|\xi|^\a, \qquad\forall\xi\in \Rn\eqdef{ellip1}$$  then its fundamental solution  $\kp$ is in $L^{\nna}(|x|\ge 1)$.\par

Note that condition \eqref{ellip1} implies the ellipticity condition \eqref{ellip0}.  As it turns out there are elliptic operators with  $p(\xi)\neq0$, $p(0)=0$ and whose fundamental solution does not satisfy the critical integrability condition. See Remark 2 after the proof of Lemma 13. 
\medskip
The proof of Theorem 11 is accomplished by showing  that the first term in \eqref{K1} is in $L^\nna(\Rn)$. From the Hausdorff-Young inequality it is enough  to prove that   $\eta/p\in L^\nsa(\Rn)$ i.e. $1/p\in L^{\nsa}(B_1)$: 

\proclaim Lemma 12.  Let $p$ be a  polynomial of even order $\alpha\in \N$ in $\R^n$, such that $p(0)=0$ and   $|p(x)|\ge c_1 |x|^\a$ for all  $x\in\R^n$, and some constant $c_1>0$. Then $p(x)^{-1}\in L^{\nsa}(B_1)$ if and only if $p$ is not homogeneous.
\par

\pf Proof of Lemma 12. Obviously, if $p$ is homogeneous of order $\a$ then $1/p$ cannot be in $L^\nsa(B_1)$. Suppose $p$ is not homogeneous and assume, without loss of generality,  that $p$ is real-valued and positive, away from $0$ (otherwise consider $|p(x)|^2$ instead of $p(x)$, and  $2\a$ instead of $\a$.)  Let $p_\a,p_{\k}$ be the highest and lowest order homogeneous parts of $p$, of orders $\alpha$ and $\k<\a$  respectively. Then we can write   $p=p_\a+q+p_\k$ and the hypotheses imply that for all $x\in \R^n$ 
$$ p_\a(x)\ge c_0|x|^\a,\qquad p_\k(x)\ge0$$ for some constant  $c_0>0$.\def\o{\omega}

Note also that $p_\kappa(\o)=0$ on a set of zero measure on $S^{n-1}$.

Write
$$\int_{|x|\le1}{1\over p(x)^\nsa} dx=\isn d\o\int_0^1{r^{n-1}\over p\big(r\o\big)^\nsa}dr=\isn d\o\int_0^{1/p_\k(\o)}{p_\k(\o)^nr^{n-1}\over p\big(rp_\k(\o)\o\big)^\nsa}dr.$$
To ease a bit the notation, for any given  $\o\in S^{n-1}$ write  $p_\a=p_\a(\o),\,p_\k=p_\k(\o)$.  We then  have
$${p_\k^nr^{n-1}\over p\big(rp_\k\o\big)^\nsa}={p_\k^n r^{n-1}\over \big(r^\a p_\k^\a p_\a+q(rp_\k\o)+r^\k p_\k^{\k+1}\big)^\nsa}={p_\k^{\nsa(\a-\k-1)} r^{n-1-{n\k\over\a}}\over \big(r^{\a-\k} p_\k^{\a-\k-1} p_\a+r^{-\k}p_\k^{-\k-1}q(rp_\k\o)+1\big)^\nsa}.
$$
We can now choose $r_0>0$ such that for all $r\le r_0$ and all $\o\in S^{n-1}$
$$r^{\a-\k} p_\k^{\a-\k-1} p_\a+r^{-\k}p_\k^{-\k-1}q(rp_\k\o)+1\ge {1\over2}$$
(recall that $\k+1\le\a$ and the lowest homogeneous part of $q$ has order greater than $\k$).
Hence we can write 
$$\int_0^{1/p_\k}{p_\k^nr^{n-1}\over p\big(rp_\k\o\big)^\nsa}dr\le2^{-\nsa}\int_0^{r_0}p_\k^{\nsa(\a-\k-1)} r^{n-1-{n\k\over\a}}dr+\int_{r_0}^{1/p_\k} {p_\k^n r^{n-1}\over \big(c_1r^\a p_\k^\a\big)^\nsa}dr
\le C(1+\big|\log p_\k\big|\big).$$
Now,  the function $\log p_\k(\o)$ is integrable on the sphere. By homogeneity it is easy to check that this is equivalent to the local integrability of $\log p_\k(x)$, which follows from this general lemma:

\proclaim Lemma 13. If $p$ is any complex-valued polynomial in $\Rn$ then the function $\log|p|$ is locally integrable in $\Rn$.\par

We have not seen this result in the literature, so we will give here  a short proof.

\pf Proof of Lemma 13. Suppose $p$ has degree $m$. By a linear transformation $x_1=x_1'$, $x_j=x_j'+\lambda_j x_1'$, $j=2,3..,n$,  we can assume that $p(x)= x_1^m+a_{m-1} x_1^{m-1}+...+a_1 x_1+a_0$, where the $a_j$ are polynomials in $x_2,...,x_n.$ If $Q=[a,b]\times Q'$ is any cube in $\Rn$, then for fixed $x_2,...,x_n\in Q'$  the polynomial $p(x)$ has $m$ complex roots $\rho_k=\rho_k(x_2,...,x_n)$, $k=1,...,m$, which are all contained inside a fixed ball of radius $R$. Then the result follows from Fubini's theorem, since
$$\int_a^b \big|\log|t-\rho|\,\big| dt\le C(a,b,R),\qquad \rho\in\C,\; |\rho|<R.$$

\endpf
\noin{\bf Remarks.}\smallskip\noin{\bf 1.}  After the first version of this manuscript was completed, Fulvio Ricci pointed out to us that Lemma 13 could also be seen as a consequence of the estimate (2.1) in [RiSt], which in our notation reads
$$\int_{|x|\le 1}|p(x)|^{-\e}dx\le A_\e \bigg(\sum_{|k|\le \a}|a_k|\bigg)^{-\e}$$
for all $\e<1/\a$. Indeed, one checks easily that the constant $A_\epsilon$ in that estimate is of type $C^\e$ and this implies the local integrability of $\log|p(x)|$. We thank Fulvio Ricci, and we also thank Peter Wagner who first pointed out to us the Proposition on p.182 of [RiSt], which contains estimate (2.1).
\smallskip\noin{\bf 2.} Lemma 12 is not valid under the weaker hypothesis $p$ elliptic, $p(0)=0$ and $p(x)\neq0$ for all $x\neq0$. An example is   the   $4-$th order elliptic  polinomial in $\R^8$
$$ p(x)=(x_8-\|x'\|^2)^2+x_8^4,\qquad x=(x_1,...,x_8)=(x',x_8)\in\R^8.$$
With a little calculation one shows that $p(x)^{-1}$ is not in $L^2(B_1)$. In fact, by Plancherel's formula this also shows that  the Fourier transform of $\eta/p$, where $\eta $ is a smooth  cutoff equal~1 on $B_1$, cannot be in $L^2(\Rn)$. This means, in view of \eqref{K1},  that the fundamental solution of the  $4$-th order operator in $\R^8$  whose symbol is  $p$, cannot be Riesz subcritical.\medskip
\smallskip
Taking into account Theorem 11 and the discussion preceding it, we can now state the following:
\smallskip
\proclaim Theorem 14. Suppose that  $P$ is any non-homogeneous elliptic differential operator of even order $\a<n$ of one the following types:
\item {i)} $P=\big(\nabla^T\A\nabla+b^T\nabla\big)^{\a/2}$, where  $\A$ is a real, symmetric and positive definite matrix and $b\in \R^n\setminus\zero$;\smallskip\item{ii)} 
 the total symbol of $P$ satisfies $p(\xi)\neq0$ for all $\xi\neq0$ and the lowest order terms of $P$ form an elliptic and homogeneous operator of order $m<\a$;
\smallskip\item{iii)} the total symbol of $P$ satisfies \eqref{ellip1}, and $\a\ge \ds{n\over2}$.
\smallskip
  Then, there exists $C$ such that for all measurable $E\subseteq\Rn$ and for all $ u\in W^{\a,\nsa}(\Rn)$ with  $\|Pu\|_\na\le 1$ we have 
$$\int_E\exp\bigg[{1\over A}|u(x)|^{\nna}\bigg]dx\le C(1+|E|)$$
with
$$A={1\over n}\isn|g_P^{}(\omega)|^\nna d\omega$$
(and with $g_P^{}$ as in \eqref{gp}).
Moreover,  for all $ u\in W^{\a,\nsa}(\Rn)$ with  $\|Pu\|_\na\le 1$ 
$$\irn\exp_{[\nsa-2]}\bigg[{1\over A}|u(x)|^{\nna}\bigg]dx\le C(1+\|u\|_\na^\na).$$
The exponential constants in both of the above inequalities are sharp.\par

\smallskip
\noin{\bf Remark.} 
In case i) the sharp exponential constant is $\gamma(\nabla^\a)(\det\A)^{\a\over 2(n-\a)}$ (cf. [FM1, Corollary 11]). This can be seen directly from the explicit formula for $K_P$ given in [OW, (2.5.3)], since $K_P(x)=c_\a(\det\A)^{-\half}|\A^{-\half}x|^{\a-n}(1+O(|\A^{-\half}x|^\sigma)$, as $x\to 0$ for some $\sigma>0$, from which it follows that   $(K_P)^*(t)=c_\a(\det\A)^{-{\a\over2n}}t^{-{n-\a\over n}}(1+O(t^\e)$), as $t\to0$, from some $\e>0$.
In cases ii)  or iii) the exponential constant can sometimes be computed explicitly, especially  when $\a=n/2$, via the Plancherel formula (cf. [FM1, formula (85)]).

\bigskip

\pf Proof.  We know that for a $P$ of the above types $\kp$ is Riesz subcritical, so that  Theorem~7 applies. 
  To prove the sharpness of the exponential constants, one would like to take the family of function $\psi_\e=\phi_\e/\|\phi_\e\|_\na$ , where the $\phi_\e$ are defined in \eqref{extremals}, and then consider $u_\e=K_P*\psi_\e$. The only problem is that one can  guarantee that $K_P*\psi_\e\in L^\na(|x|\ge 1)$  only when $\a<n/2$. For higher values of $\alpha$ one needs to first normalize the $\phi_\e$ in order to have enough vanishing moments. This is accomplished in [FM2, sect 6].
\endpf  
  \vskip 2em \eject
\centerline{\scaps 5.  Moser-Trudinger inequalities in hyperbolic space}\bigskip

In this section  we obtain the sharp Moser-Trudinger inequalities for the higher order gradients on the hyperbolic space $\Hn$, as a  consequence of Theorem 1. Below, $\Hn$ will denote the hyperbolic space 
modeled by the forward sheet of the hyperboloid
$$\Hn=\{(x_0,x_1,...,x_n)\in\R^{n+1}:\, x_0^2-x_1^2-....-x_n^2=1,\; x_0>0\},$$
 endowed with the metric induced by the form
$$[x,y]=x_0y_0-x_1y_1-...-x_ny_n$$
and with distance function $d(x,y)={\rm {arccosh}} [x,y]$. One can introduce polar coordinates on $\Hn$
via 
$$x=(\cosh r,\sinh r\, \xi),\qquad r\ge0,\,\xi\in S^{n-1}$$
and in these coordinates the metric and the volume element are  written as 
$$ds^2=dr^2+\sinh^2 r d\xi^2,\qquad d\nu(x)=(\sinh r)^{n-1}drd\xi.$$
The Laplace-Beltrami operator on $\Hn$ is denoted as $\Delta_{\sHn}$, and in polar coordinates is written  as 
$$\Delta_{\ssHn}={\p^2\over \p r^2}+(n-1)\coth r{\p\over\p r}+{1\over \sinh^2 r}\Delta_{S^{n-1}},$$
whereas the gradient $\nabla_\ssHn$ is given by 
$$\nabla_\ssHn={\partial\over\partial r}+{1\over \sinh^2 r}\nabla_{S^{n-1}}.$$

The Sobolev space $W^{\a,p}(\Hn)$ of integer order $\alpha$ is  defined in the standard way via the covariant derivatives $\nabla^k$: it is the closure of the space of $C^\infty$ functions $\phi$ such that
$$\|\phi\|_{\a,p}:=\sum_{k=0}^\alpha \|\nabla^k\phi\|_p<\infty$$
where $\|\cdot\|_p$ denotes the norm in $L^p(\Hn,\nu)$.
As it turns out, on $\Hn$ it is enough to use the highest order derivatives in order to characterize the Sobolev space. In particular, if we define the higher order gradient on $\Hn$ as
$$\nabla_{\ssHn}^\a=\cases{\nabla_{\ssHn}(-\Delta_{\ssHn})^{\a-1\over2} & if $\alpha$ odd\cr (-\Delta_{\ssHn})^{\alpha\over2} & if $\alpha$  even,\cr}$$
then one has that $\|\nabla_\ssHn^\alpha u\|_p$ is an equivalent norm on $W^{\alpha,p}(\Hn)$. In particular, note that we have the Poincar\'e-Sobolev inequality 
$$\|u\|_p\le C\|\nabla_\ssHn^\alpha u\|_p,\qquad u\in W^{\a,p}(\Hn).$$
This inequality is proved in [Mancini-Sandeep-Tintarev] in the case of the gradient in the ball model (really a consequence of Hardy's inequality) and for even $\alpha$ in [Tat].

In this setup  sharp versions of the Moser-Trudinger inequality for $W^{\alpha,\nsa}(\Hn)$ are only known in the case $\alpha=1$ for the gradient ([MS], [MST], [LT1], [LT2]), and with the same sharp constant as in the Euclidean case. In the following theorem we give the general version of this result for arbitrary $\alpha:$
\medskip\
\proclaim Theorem 15. For any integer $\alpha$  with  $0<\alpha<n$ there exists a constant $C=C(\alpha,n)$ such that for every  $u\in W^{\a,\nsa}(\Hn)$ with 
$\|\nabla_\ssHn^\alpha u\|_\na\le 1$, and for all measurable $E$ with $0<\nu(E)<\infty$ we have 
$$\int_E\exp\Big[ \gamma_{n,\a}|u(x)|^{n\over n-\a}\Big]d\nu(x)\le C\big((1+\nu(E)\big),\eqdef {Z2}$$
and
$$\int_\sHn\exp_{[\nsa-2]}\Big[ \gamma_{n,\a}|u(x)|^{n\over n-\a}\Big]d\nu(x)\le C,\eqdef {Z2'}$$
  and the constant $\gamma_{n,\a}$ is sharp.

\par

\pf Proof. If $\alpha$ is even, the operator $(-\Delta_\ssHn)^{\a\over2}$ has a fundamental solution given by a  kernel of type $H_\a\big(d(x,y)\big)$, where $H_\a$ is positive and satisfies
$$H_\a(\rho)=c_\a \rho^{\a-n}+O(\rho^{\a-n+\e}),\qquad  \rho<1\eqdef{Z3}$$
(with the same $c_\a$ as in the Euclidean Riesz potential), and  $$H_\a(\rho)\le c_\a'\rho^{-1+{\a\over2}}e^{-(n-1)\rho}, \qquad \rho\ge1,\eqdef{Z4}$$
some $c_\a'>0$. These asymptotic estimates follow in a straightforward manner from the known formula for the fundamental solution of the Laplacian (see for example [CK])
$$H_2(\rho)={1\over\omega_{n-1}} \int_\rho^\infty {dr\over (\sinh r)^{n-1}}$$
using iterated integrations and the known addition formulas for the Riesz potential on $\R^n$. (In [BGS]  asymptotic formulas are derived  for general $\a$, using the  Fourier transform.)

It is now easy to check that \eqref{Z3} implies that in the measure space $(\Hn,\nu)$ we have  $$H_\a^*(t)=c_\a|B_1|^{n-\a\over n}t^{-{n-\a\over n}}+O(t^{-{n-\a\over n}+\e}),\quad t\le 1$$
while \eqref{Z4} implies that $H_\a\big(d(\cdot,O)\big)\in L^{\nna}\cap L^\infty\big(\{x:d(x,O)\ge1\}, \nu\big)$ (where \break$O=(1,0,...,0)$) and hence
$$\int_1^\infty (H_\a^*)^{n\over n-\a}dt<\infty.$$
Thus, we are in a position to apply Theorem 7 in  order to obtain \eqref{Z2} for $\alpha$ even, simply by writing $u(x)=\int H_\a(d(x,y))f(y)d\nu(y)$, with $f= \nabla_\ssHn^\a u$, for any $u\in C_0^\infty(\Hn)$.

If $\alpha$ is an odd integer, then we write 
$$u(x)=\int_\sHn \nabla_\ssHn H_{\a+1}(d(x,y))\cdot f(y)d\nu(y),\qquad f=\nabla_\ssHn^\a u$$
and use asymptotic estimates for $|\nabla_\ssHn H_{\a+1}|$,  which turn out to be the same exact estimates as in \eqref{Z3}, \eqref{Z4}, with $(n-\a-1)c_{\a+1}$ instead of $c_\a$.

The proof of the sharpness statement is identical to the one in the Euclidean case, namely we  let $v_\e$   to be  a   smoothing of the radial function
$$\cases{0 & if $r\ge{3\over4}$\cr \log{1\over r} & if $2\e\le r\le {1\over2}$\cr
\log{1\over\e} & if $r\le \e$.\cr}$$
Using local calculations as in [F, Prop. 3.6]  it is a routine task to check that 
if $\alpha$ is even then
 $$ \|\nabla_\ssHn^\a v_\e\|_\na^\na=\omega_{n-1}^{-{n-\a\over \a}}c_{\a}^{-\nsa}\log{1\over\e}+O(1),  $$
whereas if $\alpha$ is odd then the same estimate holds with $(n-\a-1)c_{\a+1}$ in place of $c_\alpha$.
From this estimate it is then clear that the exponential integral evaluated at the functions $u_\e=v_\e/\|\nabla_\ssHn^\a v_\e\|_\na$ can be made arbitrarily large  if the exponential constant is larger than $ \gamma_{n,\a}$. Note also that $\|v_\e\|_\na\le C$, so that $\|u_\e\|_\na\to0$ with $\e$, and the sharpness statement for the regularized inequality on $\Hn$ follows as well. \endpf\vskip1em

\bigskip
\centerline{\scaps 6. Connections with  the Poincar\'e inequality:}
\centerline{\scaps The Moser-Trudinger inequality on Agmon-Souplet domains}
\bigskip
In this section we are concerned with the validity of the Moser-Trudinger inequality 
$$\int_{\Omega}\exp_{[\nsa-2]}\Big[\gamma(\nabla^\a)|u(x)|^\nna\Big]dx\le C,\qquad \forall u\in W_0^{\a,\na}(\Omega) {\hbox{ such that }}\; \|\nabla^\a u\|_\na\le 1\eqdef{MTS}$$
where $\a$ is an integer in $(0,n)$, and $\Omega$ an open set in $\Rn$. From \eqref{MTS0} of Corollary 10 we know that if $\Omega$ is Riesz subcritical and the Poincar\`e inequality holds in the form 
$$\|u\|_\na\le C\|\nabla^\a u\|_\na,\qquad \forall u\in W_0^{\a,\na}(\Omega)\eqdef{P1}$$
then \eqref{MTS} also holds.
On the other hand,  for $\alpha=1$ Battaglia and Mancini [BM] proved that \eqref{MTS} holds if and only if  \eqref{P1} holds.

One direction of this result is in some sense an artifact of the exponential regularization. Indeed it is clear from  \eqref{exp2} of Lemma 4 that if \eqref{MTS} holds under the hypothesis $\|\nabla^\a u\|_\na\le~1$, then $\|u\|_\na\le C$ provided that $\na$ is an integer.
The interesting part is the  reverse implication: when is it true that the Poincar\'e inequality for $p=\na$ implies the Moser-Trudinger inequality \eqref{MTS}? We know this fact for $\alpha=1$, by the above-mentioned result in [BM], and for general $\alpha$ in the case of Riesz subcritical domains. The question remains open for the general case $\alpha>1$, since there are are indeed domains $\Omega$ which are  not Riesz subcritical and on which  the Poincar\`e inequality holds for any $p>1$:
$$\|u\|_p\le C\|\nabla^\a u\|_p,\qquad u\in W_0^{\a,p}(\Omega)\eqdef{Poincare}.$$
\noin To describe such domains in some generality we recall a few definitions.


The {\it  inradius} of a domain (nonemtpy, open) $\Omega\subseteq \Rn$  is defined as 
$$\rho(\Omega)=\sup \{r>0: \;\Omega{\hbox{ contains a ball of radius }}r\}.$$
whereas the strict inradius is defined as
$$\rho'(\Omega)=\sup \{R>0: \;\forall\e>0,\,\exists B(a,R) {\hbox{ such that }} B(a,R)\cap \Omega^c {\hbox{ contains no ball of radius }}\e\}.\eqdef{si1}$$
Note that the set appearing in the  definition of $\rho'(\Omega)$ is an open interval of type $(0,a)$, so that 
$$\eqalign{\rho'(\Omega)=\inf \{R>0: \;\exists\e>0&  { \hbox{  such that for any  ball }}  B(a,R),\cr &\hskip 1.8em B(a,R)\cap \Omega^c {\hbox{ contains a ball of radius }}\e\}.\cr}\eqdef{si2}$$
The notion of strict inradius as stated is due to Souplet ([So], [QS]), who proved that a {\it sufficient} condition for the validity of the Poincar\`e inequality 
\eqref{Poincare}
 is that $\rho'(\Omega)<\infty$ (Souplet proved it for $\alpha=1$, but from there it is  easy to  extend it  to any $\a\in\N$). This result was due to Agmon in [Ag] in the case $p=2,\a=1$ under the condition
$$\exists R_1,\e_1 {\hbox { such that }} \forall x\in\Omega,\,\exists x^*\in\Omega^c:\quad |x-x^*|\le R_1 \;{\hbox { and }}\; B(x^*, \e_1)\subseteq \Omega^c\eqdef{si3}$$
which is  equivalent to $\rho'(\Omega)<\infty$ on $\Rn$. In this paper we will call $\Omega$ an {\it Agmon-Souplet domain} if $\rho'(\Omega)<\infty$, i.e.  if $\Omega$ satisfies condition \eqref{si3}.

Examples of Agmon-Souplet domains include bounded smooth domains (in which case $\rho(\Omega)=\rho'(\Omega)$), domains contained in a ``strip'', complements of periodic nets of balls whose radii is bounded below by a positive number. It is worth observing that in such domains we cannot expect the Riesz potential to be continuous, nor the Adams inequality to hold. For example, if $\Omega=\Rn\setminus \bigcup_{\m\in\sZ^n}B(\m,\e_0)$, some fixed small $\e_0>0$, then the Agmon-Souplet condition is verified, the Poincare' inequality holds, however, in the notation of section 3, $\G(x,r)\ge C r^n$ for any fixed $x\in\Omega$ and for $r\ge r_0>0$, hence  condition  \eqref{G} of Theorem~9  is certainly not met, hence the set is not Riesz subcritical and Adams' inequality fails.
However, we are able to show, as a nice application of Theorem 1, that in domains satisfying the Agmon-Souplet condition the Moser-Trudinger inequality actually holds:

\medskip
\proclaim Theorem 16. If $\Omega$ is a domain in $\Rn$ such that $\rho'(\Omega)<\infty$, then \eqref{MTS} holds, 
and the exponential constant is sharp.\par

\pf Proof.  Let us consider first the case $\a$ even. If $u\in C_c^\infty(\Omega)$ with $\|\nabla u\|_\na\le 1$ then $u= c_\a I_\a f$, where $f=\nabla^\a u=(-\Delta)^{\a/2} u$, compactly supported inside $\Omega$. The point is that it is possible to normalize the Riesz kernel so that the critical integrability condition at infinity  is satisfied, without interfering with the local asymptotics. Indeed,  we can cover $\Omega$ with countably many balls $B(x_j,R_1)$, with $x_j\in \Omega$. If  for each $j$ we pick $x_j^*\in \Omega^c$ according to  \eqref{si3}, then  $B(x_j,R_1)\subseteq B(x_j^*,2R_1)$ and   $B(x^*_j,\e_1)\subseteq \Omega^c$. Since $I_\a f(x_j^*)=0$ for all $j$, then for each $j$ we can write
$$I_\a f(x)=\int_\Omega \bigg(|x-y|^{\a-n}-|x^*_j-y|^{\a-n}\bigg)f(y)dy,\qquad \forall x\in B(x_j^*,2R_1)\cap \Omega.$$
Let
$$B_1=B(x_1^*,2R_1),\quad B_j=B(x_j^*,2R_1)\setminus\bigcup_{k=1}^{j-1}B(x_k^*,2R_1), \;\;j\ge2$$
then $\bigcup_j B_j\cap \Omega=\Omega$ and the $B_j$ are disjoint. For each $(x,y)\in \Omega\times\Omega$ define 
$$K(x,y)=c_\a\sum_{j=1}^\infty \Big(|x-y|^{\a-n}-|x^*_j-y|^{\a-n}\Big)\chi_{\Omega\cap B_j}^{}(x)=c_\a|x-y|^{\a-n}-c_\a\sum_{j=1}^\infty |x^*_j-y|^{\a-n}\chi_{\Omega\cap B_j}^{}(x) $$
then $u(x)=c_\a I_\a f(x)=\int_{\Omega}K(x,y)f(y)dy$ for each $x\in \supp f$, so it is enough to show that $K_1^*$ and $K_2^*$ satisfy conditions \eqref{106}, \eqref{107}, \eqref{108} of Theorem 1.

If $|x-y|\le \e_1$ then $|K(x,y)|\le c_\a|x-y|^{\a-n}$. 
On the other hand, for $x\in B_j\cap \Omega$ and $y\in \Omega$ we have
$$|K(x,y)|=c_\a\big||x-y|^{\a-n}-|x_j^*-y|^{\a-n}\big|\le c_\a2R_1(n-\a)\max\big\{|x-y|^{\a-n-1},|x_j^*-y|^{\a-n-1}\big\}$$
from which we deduce that if $|x-y|\ge \e_1$, then $|K(x,y)|\le C\e_1^{\a-n-1}R_1$, and if $|x-y|\ge 2R_1$ then
$|x_j^*-y|\ge \half |x-y|$ and $|K(x,y)|\le CR_1|x-y|^{\a-n-1}$.
From these estimates it is straightforward to check that $K_1^*$ and $K_2^*$ satisfy \eqref{106}, \eqref{107}, and \eqref{108} of Theorem 1, with  $\beta=\nna$, $A^{-1}=\gamma(\nabla^\a)$, $\sigma=1$. 
\smallskip
The proof  in the case $\a$ odd is similar, starting from the identity
$$u(x)=J_\a f(x):= c_{\a+1}(n-\a-1)\int_\Omega|x-y|^{\a-n-1}(x-y)\cdot f(y)dy,\quad f=\nabla(-\Delta)^{\a-1\over2}u.$$
In this case we normalize the kernel of $J_\a$ by letting, for $(x,y)\in \Omega\times\Omega$
$$K(x,y)=c_{\a+1}(n-\a-1)\sum_{j=1}^\infty \Big(|x-y|^{\a-n-1}(x-y)-|x_j^*-y|^{\a-n-1}(x_j^*-y)\Big)\chi_{ \Omega\cap B_j}^{}(x).$$
If $|x-y|\le \e_1$ then use $|K(x,y)|\le c_{\a+1}(n-\a-1)(|x-y|^{\a-n}+\e_1^{\a-n})$, whereas if $|x-y|\ge\e_1$ then use for $x\in \Omega\cap B_j$, $y\in \Omega$
$$\eqalign{|K(x,y)|\le &C\Big||x-y|^{\a-n-1}(x-y)-|x-y|^{\a-n-1}(x_j^*-y)\Big|+\cr&\hskip12em +C\Big||x-y|^{\a-n-1}(x_j^*-y)-|x_j^*-y|^{\a-n-1}(x_j^*-y)\Big|\cr&\le CR_1|x-y|^{\a-n-1}+C\e_1\Big||x-y|^{\a-n-1}-|x_j^*-y|^{\a-n-1}\Big|\cr& \le C(R_1+\e_1)\max\big\{|x-y|^{\a-n-1},|x_j^*-y|^{\a-n-1}\big\}\cr}$$
and the same estimates as in the case $\a$ even apply.

The proof of the sharpness statement is the same one as in the classical case of bounded domains.\endpf\bigskip

 Note that it is possible to construct  domains for which the Poincar\'e inequality \eqref{Poincare} holds  for all $p\ge 1$ and which are not satisfying the Agmon-Souplet condition. Here's an outline of this construction. 
First, we note the following variation of Theorem 9:

\medskip
\proclaim Proposition 17. If $0<\a<n$, and  $\Omega\subseteq\Rn$ is measurable and such that (using the notation in \eqref{lambda})
$$\sup_{x\in \Omega}\int_0^\infty {\G(x,r)\over r^{n+1-\a}}dr<\infty,\eqdef{GG}$$
then the Riesz potential is continuous from $L^p(\Omega)$ to $L^p(\Omega)$. In particular, if $\Omega$ is open, and \eqref{GG} holds for at least  $\alpha=1$, \eqref{Poincare} holds for any integer $\alpha$ and any $p\ge1$.\par
\pf Proof. By a standard result the continuity follows from 
$$\sup_{x\in\Omega}\int_\Omega |x-y|^{\a-n}dy<\infty,\eqdef{Young}$$
which is the same as \eqref{GG}, since
$$\int_\Omega |x-y|^{\a-n}dy=\int_0^\infty|\{y\in\Omega: |y-x|^{\a-n}>t\}|dt.$$\endpf

Condition \eqref{GG} is meaningful only for large $r$, say $r\ge 1$, (since the integrand is integrable around 0), and it is clearly stronger than \eqref{G}.  For example,  when $\Omega=(a,b)\times \R^{n-1}$ then \eqref{G} holds but \eqref{GG} holds  only  for $0<\a<1$.

Now consider the open subset  $A\subseteq (0,\infty)$ obtained by removing from each $(k,k+1)$  $2^k-1$ equally spaced intervals  of length $\d_k<2^{-k}$. Specifically, let 
$$I_k^j=\Big[k+{j\over 2^{k}}-{\d_k\over2},k+{j\over 2^{k}}+{\d_k\over2}\Big],\quad k\in\N,\; j=1,...,2^k-1$$
$$A=(0,\infty)\setminus \bigcup_{k=1}^\infty\bigcup_{j=1}^{2^k-1}I_k^j.$$
If we let  $S_m=\sum_0^m \big(1-(2^k-1)\delta_k\big)$, and assuming that $(2^k-1)\delta_k$ is increasing,  then it is easy to see that in $\R$
$$S_{[r]-1}\le \Lambda_A(0,r)\le S_{[r]},\qquad \Lambda_A(x,r)\le S_{2[r]+1}\le \Lambda_A(0,2r+2),\qquad r\ge 1, x>0.$$

Now condition \eqref{GG} holds if and only if it holds with $\G(x,r)$ replaced by 
$$\wtilde \G(x,r)= |\Omega\cap Q(x,r)|$$
where $Q(x,r)$ is the cube centered at $x$ and with sidelength $2r$. If we set
$$\Omega=A\times A\times ....\times A=\prod_{j=1}^n A\subseteq \Rn$$
then $\Omega$ is an open set in $\Rn$ for which the Agmon-Souplet condition is not satisfied (for each $R,\e>0$ a point  $x$ can be  chosen far enough inside $\Omega$ so that  $B(x,R)\cap\Omega^c$ is not empty, but it does not contain any ball of radius $\e$). Yet, Poincar\'e's inequality holds in such domain,  since
$$\wtilde \G(x,r)\le \big(S_{2[r]+1}\big)^n,\qquad \wtilde \G(0,r)\ge \big(S_{[r]-1}\big)^n$$ and
 one can choose $\delta_k$ satisfying the above conditions in a way that \eqref{GG} is satisfied for all $\alpha<n$, with $|\Omega|$ either  finite or infinite (for example by choosing $\d_k$ so that $(1-(2^k-1)\d_k)=Ak^{-\gamma}$, for $\gamma>0$ and  suitable $A>0$).
\medskip




\vskip2em

\centerline{\scaps 6. Appendix}
\bigskip
\centerline{\bf Proof of Lemma 5 (Improved O'Neil's Lemma)}
\bigskip

The given hypothesis \eqref{on1} are equivalent to the weak-type estimates 
$$\sup_{x\in N} \lambda_1(x,s)\le D s^{-\beta},\qquad \sup_{y\in M} \lambda_2(y,s)\le Bs^{-\sigma\beta},$$
 and a result due to Adams [A2]  gives
$$ s\,\nu\big(\{x: Tf(x)>s\}\big)^{1\over q}\le {q^2\over \sigma\beta (q-p)}D^{1-{1\over p}}B^{1\over q} \|f\|_p\eqdef{(21)}$$
or
$$ (Tf)^*(t)\le C t^{-{1\over q}} \|f\|_p,\qquad \forall t>0.\eqdef{(22)}$$
In particular this means that $Tf$  is well defined on $L^p(M)$ and  bounded from $L^p(M)$ to $L^{q,\infty}(N)$.

\proclaim Claim. If $\mu(\supp f)=z$ and $0\le f(z)\le \alpha$, and if $k_1^*,k_2^*$ satisfy conditions \eqref{on1}, then 
$$Tf(x)\le\alpha \int_0^z k_1^*(x,v)dv,\qquad  \nu-{\hbox{a.e. }} x\in N \eqdef{(23)}$$
$$(Tf)^{**}(t)\le C\, \alpha\, z^{{1\over p}} t^{-{1\over q}},\qquad \forall t>0.\eqdef{(24)}$$\par
\noin Estimate \eqref{(23)} is an improvement of the corresponding estimate given in [FM1, Lemma~2], which was given as $(Tf)^{**}(t)\le \a z k_1^{**} (z)$. \smallskip
Assuming the Claim, the proof of the lemma proceeds as follows.
For fixed $t,\tau>0$, pick  $\{y_n\}_{-\infty}^\infty$ such that $y_0=f^*(\tau),\,y_n\le y_{n+1}, \,y_n\to+\infty$
as $n\to+\infty$, and $y_n\to0$ as $n\to-\infty$.  Then 
$$f(y)=\sum_{-\infty}^\infty f_n(y)\quad{\hbox{where}}\quad f_n(y)=\cases{0 &if $\;f(y)\le y_{n-1}$\cr f(y)-y_{n-1} 
&if $\;y_{n-1}<f(y)\le y_n$ \cr y_n-y_{n-1} & if $\;y_n<f(y).$\cr}$$

Observe that supp$f_n\subseteq E_n:=\big\{y: f(y)>y_{n-1}\big\}$,  $\;\mu(E_n)=\lambda_f(y_{n-1})$, with 
$$\lambda_f(s):=\mu(\{y\in M:\; |f(y)|>s\}),$$
and also  $\;0\le f_n(y)\le y_n-y_{n-1}$. Write 
$$f=\sum_{-\infty}^0f_n+\sum_1^{\infty}f_n=g_1+g_2.$$
By the  subadditivity of $(\cdot)^{**}$ ([BS], Ch. 2, Thm. 3.4), and  using  \eqref{(24)} we obtain
$$(Tg_2)^{**}(t)\le \sum_1^\infty (Tf_n)^{**}(t)\le Ct^{-{1\over q}}\,\sum_1^{\infty}(y_n-y_{n-1})\big(\lambda_f(y_{n-1})\big)^{{1\over p}}$$
so that taking the inf over all such $\{y_n\}$ we get
$$\eqalign{&(Tg_2)^{**}(t)\le Ct^{-{1\over q}}\int_{f^*(\tau)}^\infty \big(\lambda_f(s)\big)^{{1\over p}} ds=- Ct^{-{1\over q}}\int_0^\tau\big(\lambda_f(f^*(u))\big)^{{1\over p}}d\,f^*(u)\cr&\le-  Ct^{-{1\over q}}\int_0^\tau u^{{1\over p}} d\,f^*(u)= Ct^{-{1\over q}}\bigg(-u^{{1\over p}}f^*(u)\Big|_0^\tau+{1\over p}\int_0^\tau u^{-1+{1\over p}}f^*(u)du\bigg)\cr& \le{ Ct^{-{1\over q}}\over p}\int_0^\tau u^{-1+{1\over p}}f^*(u)du. \cr}  $$
(The last inequality follows since $f\in L^{\beta'}\Longrightarrow t^{{1\over \beta'}}f^*(t)\to0$, as $t\to0$.)

Using  \eqref{(23)}, for $\nu-$a.e. $x\in N$ we have
$$Tg_1(x)\le \sum_{-\infty}^0 Tg_1(x)\le \sum_{-\infty}^0 (y_n-y_{n-1})\int_0^{\lambda_f(y_{n-1})}
k_1^*(x,u)du$$
and so, arguing as above, 
$$\eqalign{&Tg_1(x)\le\int_0^{f^*(\tau)}ds\int_0^{\lambda_f(s)}k_1^{*}(x,v)dv=
-\int_\tau^\infty df^*(u)\int_0^{\lambda_f(f^*(u))}k_1^{*}(x,v)dv\cr&\le-\int_\tau^\infty df^*(u)\int_0^{u}k_1^{*}(x,v)dv=
-f^*(u) \int_0^u k_1^{*}(x,v)dv \bigg|_\tau^\infty+\int_\tau^\infty k_1^*(x,u)f^*(u)du\cr&\le f^*(\tau)\int_0^\tau k_1^{*}(x,v)dv+\int_\tau^\infty k_1^*(x,u)f^*(u)du\cr& \le \tau^{1-{1\over p}}\int_0^\tau Dv^{-{1\over\beta}}dv\int_0^\tau f^*(u)u^{-1+{1\over p}}du+\int_\tau^\infty k_1^*(x,u)f^*(u) du\cr&
\le C\,\tau^{-{\sigma\over q}}\int_0^\tau f^*(u)u^{-1+{1\over p}}du+\int_\tau^\infty k_1^*(x,u)f^*(u)du
}.\eqdef{25b}$$
Finally,
$$Tf^{**}(t) \le (Tg_1)^{**}(t)+(Tg_2)^{**}(t)\le \|Tg_1\|_\infty +Ct^{-{1\over q}} \int_0^\tau  u^{-1+{1\over p}} f^*(u)du$$
which, together with \eqref{25b}, implies (19). 



\medskip
\smallskip\noin{\bf Proof of Claim.}  Estimate \eqref{(24)} is an immediate consequence of the weak-type estimate \eqref{(22)}. To show \eqref{(23)}, let $r>0$ and   set 
$$k_r(x,y)=\cases{k(x,y) & if $\;k(x,y)\le r$\cr\cr r & otherwise,\cr}\qquad k(x,y)=k_r(x,y)+k^r(x,y).$$
so that 
$$Tf(x)=\im k_r(x,y)f(y)d\mu(y)+\im k^r(x,y )f(y)d\mu(y)=h_1(x)+h_2(x).$$
Then,
for every given $x$
$$ h_2(x)\le \|f\|_\infty^{}\im k^r(x,y)d\mu(y)\le \alpha \int_r^\infty \lambda_1(x,s)ds,\eqdef{(26)}$$
$$ h_1(x)\le \|f\|_1^{}\sup_y k_r(x,y)\le \alpha z r,\eqdef{(27)}$$
so that letting $r=k_1^*(x,z)$ in \eqref{(26)} and \eqref{(27)} leads to

$$\eqalign{Tf(x)\le
\alpha z \,k_1^*(x,z)+\alpha\int_{k_1^*(x,z)}^\infty \lambda_1(x,s)ds=\alpha\int_0^z k_1^*(x,u)du\cr}$$
which is \eqref{(23)}.

 \endpf
\centerline{\bf Proof of inequality \eqref{AG}.}

\bigskip
Let 
$$L(\eta)=\bigg(\int_\eta^\infty \phi(\xi)^{\b'}d\xi\bigg)^{\bpi}\le \|\phi\|_{\b'}\le1.$$
In what follows we will repeatedly make use of the following inequalities
$$(a+b)^\b\le a^\b+\beta 2^{\b-1}(a^{\b-1} b+b^\b),\qquad ab\le {a^\b\over\b}+{b^{\b'}\over \b'},\qquad a,b\ge0$$
$$\bigg(\sum_1^m a_k\bigg)^\b\le m^\b\sum_1^m a_k^\b,\qquad  a^{\bpi}\le 1+a.$$
Note that if $0\le z_1\le z_2 \le \eta$ we have
$$\sup_{x\in N} \int_{z_1}^{z_2} g(x,\xi,\eta)^\b d\xi\le\int_{z_1}^{z_2} (1+ H(1+|\xi-\log \ts|)^{-\gamma}\big)^\b d\xi\le z_2-z_1+ C_3.$$
Also, for $z\ge \eta\ge0 $ we have
$$\sup_{x\in N}\int_z^\infty g(x,\xi,\eta)^\b d\xi\le \int_\eta^\infty  C_2^\beta e^{{\b\over q}(\eta-\xi)}d\xi={q\over \b}\, C_2^\beta=C_4$$
and 
$$\sup_{x\in N} \int_{-\infty}^{0} g(x,\xi,\eta)^\b d\xi\le {\s  } J_\t.$$
Next, we note that for $\eta>0$ H\"older's inequality implies
$$\eqalign{F(\eta)\ge \eta-\sup_{x\in N}\int_{-\infty}^\infty g(x,\xi,\eta)^\b d\xi&\ge \eta-\Big(\eta+ C_3+C_4+{\s  } J_\t\Big)\cr&\ge- C_3-C_4-{\s  } J_\t.\cr}\eqdef{2}$$
From now on let
$$d^*=C_3+C_4+{\s} J\t.$$
\def\imp{\Longrightarrow}
Now let for $\lambda\in\R$
$$E_\lambda=\{\eta\ge 0 : F(\eta)\le \lambda\}$$
and let us prove that there exists $C_5$ such that 
$$|E_\lambda|\le C_5 (|\lambda|+d^*)\eqdef{3a}$$
Proceeding as in [A1] and [FM1], it is enough to prove that there exists $C_6$ such that 
for any $\lambda\in\R$ 
$$\eta ,\eta'  \in E_\lambda,\quad \eta'  >\eta >C_6(|\lambda|+d^*)\imp \eta'  -\eta \le C_6(|\lambda|+d^*)\eqdef{3}$$
indeed, if that is the case, then
$$\eqalign{|E_\lambda|&=\big|E_\lambda\cap\{\eta: \eta\le C_6(|\lambda|+d^*)\}\big|+\big|E_\lambda\cap\{\eta: \eta> C_6(|\lambda|+d^*)\}\big|\cr& \le  C_6(|\lambda|+d^*)+\sup_{\eta'  >\eta >C_6(|\lambda|+d^*)\atop \eta ,\eta'  \in E_\lambda} (\eta'  -\eta )\cr&\le 3C_6(|\lambda|+d^*)\cr} $$
which implies \eqref{3a}.

We now prove \eqref3. If $\eta ,\eta'  \in E_\lambda$ and $|\lambda|<\eta <\eta'  $, then $F(\eta'  )\le\lambda$, so that 

$$\eqalign{(&\eta'  -\lambda)^\bi\le\sup_{x\in N}\int_{-\infty}^\infty g(x,\xi,\eta')\phi(\xi)d\xi=\sup_{x\in N}\bigg(\int_{-\infty}^{\eta }+\int_{\eta }^{\eta'  }+\int_{\eta'  }^\infty\bigg)\cr&\le
\sup_{x\in N}\bigg(\int_{-\infty}^{\eta }\!g(x,\xi,\eta')^\b d\xi\bigg)^{\bi}\!\!+
\bigg[\sup_{x\in N}\bigg(\int_{\eta }^{\eta'  } \!g(x,\xi,\eta')^\b d\xi\bigg)^{\bi}\!\!+\sup_{x\in N}\bigg(\int_{\eta'  }^\infty\! g(x,\xi,\eta')^\b d\xi\bigg)^{\bi}\bigg] L(\eta )\cr&\le
(\eta +d^*)^\bi +\big[(\eta'  -\eta +C_3)^{\bi}+C_4^\bi\big]L(\eta )
\cr}$$
(note that $g(x,\xi,\eta)=g(x,\xi,\eta')$ if $\xi\le\eta\le \eta'$)
from which we deduce
$$\eqalign{&\eta'  -\lambda\le \eta +d^*+\b 2^{\b-1}\big[(\eta +d^*)^\bpi\big((\eta'  -eta_1+C_3)^\bi+C_4^\bi\big)L(\eta )+\cr&\hskip22em + \big((\eta'  -\eta +C_3)^\bi+C_4^\bi\big)^\b L(\eta )^\b\Big]\cr&\le \eta +d^*+\b 2^{\b-1}\big[(\eta +d^*)^\bpi\big((\eta'  -\eta +C_3)^\bi+C_4^\bi\big) L(\eta )+2^\b (\eta'  -\eta +C_3)L(\eta )^\b+2^\b C_4 \Big].
\cr}\eqdef 4$$
Now we show that there exists $C_7$ such that 
$$(\eta +d^*)L(\eta )^{\b'}\le C_7\big(|\lambda|+d^*\big)\eqdef 5$$
Indeed, proceeding as above\def\bp{{\b'}}
$$\eqalign{\eta -\lambda&\le\bigg(\sup_{x\in N}\int_{-\infty}^\infty g(x,\xi,\eta )\phi(\xi)d\xi\bigg)^\b=\sup_{x\in N}\bigg(\int_{-\infty}^{\eta }+\int_{\eta }^\infty\bigg)^\b\cr& \le\Big(
(\eta +d^*)^\bi\big (1-L(\eta )^{\b'}\big)^\bpi+C_4^\bpi L(\eta )\Big)^\b
\cr
&\le (\eta +d^*)\big (1-L(\eta )^{\b'}\big)^{\b\over\bp}+\b 2^{\b-1}\big[(\eta +d^*)^\bpi\big (1-L(\eta)^{\b'}\big)^{\b-1\over\bp}C_4^\bi L(\eta )+C_4 L(\eta )^\b\big]\cr&
\le(\eta +d^*)\Big (1-\min\Big\{1,{\b\over\bp}\Big\}L(\eta )^{\b'}\Big)+C_4^\bi\b 2^{\b-1}(\eta +d^*)^\bpi L(\eta )+\b 2^{\b-1}C_4
\cr}$$
or 
$$-\lambda\le d^*-\min\Big\{1,{\b\over\bp}\Big\}(\eta +d^*)L(\eta )^\bp+C_4^\bi \b2^{\b-1}(\eta +d^*)^\bpi L(\eta )+\b2^{\b-1}C_4.$$
Letting $z=(\eta +d^*)^\bpi L(\eta )$ the last inequality can be written as
$$ z^\bp\le C_8(z+ \lambda+d^*)\le {C_8^\b\over\b}+{ z^{\b'}\over\b'}+C_8(|\lambda|+d^*)$$
which proves \eqref5. Back to \eqref 4 
$$\eqalign{\eta'  -\eta &\le \lambda+d^*+(\eta'  -\eta +C_3)^\bi\big[\b 2^{\b-1}(\eta +d^*)^\bpi L(\eta )\big]+\b2^{\b-1}C_4^\bi(\eta +d^*)^\bpi L(\eta )\cr&\hskip17em+\b2^{2\b-1}(\eta'  -\eta +C_3)L(\eta )^\b+C_4 \b2^{2\b-1}\cr
&\le \lambda+ d^*+{\eta'  -\eta +C_3\over\b}+{\big(\b2^{\b-1}\big)^\bp(\eta +d^*)L(\eta )^\bp\over\b'}+\b2^{\b-1}C_4^\bi C_7^\bpi(|\lambda|+d^*)^\bpi\cr&\hskip17em +\b2^{2\b-1}(\eta'  -\eta +C_3)L(\eta )^\b+C_4 \b2^{2\b-1}
\cr&\le {\eta'  -\eta \over \b}+C_9(|\lambda|+d^*)+C_{10}(\eta'  -\eta )L(\eta )^\b
\cr
}$$
so that
$${\eta'  -\eta \over\bp}\le C_9(|\lambda|+d^*)+C_{11}(\eta'  -\eta )\bigg({|\lambda|+d^*\over \eta +d^*}\bigg)^{\b\over\bp}. $$
Taking $\eta >C_6(|\lambda|+d^*):=(2\b' C_{11})^{\bp\over\b}(|\lambda|+d^*)$  gives 
$\eta'  -\eta \le 2\b'C_{9}\big(|\lambda|+d^*\big), $ which is \eqref 3.

\smallskip
To complete the proof we now estimate
$$\eqalign{\int_0^\infty &e^{-F(\eta)}d\eta=\int_{-d^*}^\infty |E_\lambda|e^{-\lambda}d\lambda\le \int_{-d^*}^\infty \big(C_5(|\lambda|+d^*)\big)e^{-\lambda}d\lambda\le C_{12}d^*e^{d^*}\le C_{13}\Big(1+{\s }J_\t\Big) e^{{\s }J\t}. 
\cr}$$
\endpf
\bigskip
\centerline{ \bf Proof of estimates \eqref{stimaG1} and  \eqref{stimaG2}}
\medskip First note that for each real $t> -\d_0$ there is a unique integer $m_t\ge0$ such that 
$$h(m_t)-\d_0<t\le h(m_t+1)-\d_0.$$
If 
$\wtilde\G(x,r)=|\Omega\cap Q(x,r)|$ where $Q(x,r)$ is the cube of center $x$ and side length $2r$, then $\wtilde \G(x,r/\sqrt n)\le \G(x,r)\le \wtilde \G(x,r)$. If  $r>\d_0$ then the definition of $m_r$ gives
 $$m_r^n|B_{\d_0}|\le \wtilde\G(0,r)\le (m_r+1)^n|B_{\d_0}|.$$
Choosing $r>h(4)-\d_0$ we also have $m_r\ge 3$ and 
$$\eqalign{2^{-n}\big(h^{-1}(r+\d_0)\big)^n|B_{\d_0}|\le \Big({m_r+1\over2}\Big)^n |B_{\d_0}|\le m_r^n|B_{\d_0}|&\le\wtilde \G(0,r)\le 2^n m_r^n|B_{\d_0}|\cr& \le 2^n \big(h^{-1}(r+\d_0)\big)^n|B_{\d_0}|\cr}$$ which gives \eqref{stimaG1}.

To prove \eqref{stimaG2} note that if  $r>\d_0$ and $t>-\d_0$ we have $t+r\le h(m_t+m_r+2)-2\d_0$ (if not, then 
$r>h(m_t+m_r+2)-t-2\d_0\ge h(m_t+m_r+2)-h(m_t+1)-\d_0\ge h(m_r+1)-\d_0$, using that $h(m+1)-h(m)$ is increasing).

Likewise, if $t>h(m_r+1)-\d_0$ then $m_t\ge m_r+1$ and 
$$t-r>h(m_t)-h(m_r+1)\ge h(m_t-m_r-1),$$
 so that  the number of integers $m$ such that the interval $\big(h(m)-\d_0,h(m)+\d_0\big)$ is inside $(t-r,t+r)$ does not exceed $(m_t+m_r+2)-(m_t-m_r-1)=2m_r+3$. The same is true if  $-\delta_0<t\le h(m_r+1)-\d_0$, (which is less than $r$) since $m_t+m_r+2\le 2m_r+2$.    It follows that the number of open cubes centered at $C_\m$ and with side length $2\d_0$ inside a cube of center $x$ and side length $2r$ does not exceed $(2m_r+3)^n$. This implies that for $r>h(4)-\d_0$ 
$$\G(x,r)\le \wtilde \G(x,r)\le |B_{\d_0}|(2m_r+3)^n\le 3^nm_r^n|B_{\d_0}|\le 3^n\wtilde \G(0,r)\le 3^n \G(0,r\sqrt n)$$
which is \eqref{stimaG2}.
\bigskip

\vskip2em 
\centerline{\bf References }
\vskip1em

\item{[A1]} Adams D.R. {\sl
A sharp inequality of J. Moser for higher order derivatives},
Ann. of Math. {\bf128} (1988), no. 2, 385--398. 
\smallskip
\item{[A2]} Adams D.R., {\sl A trace inequality for generalized potentials},
Studia Math. {\bf 48} (1973), 99--105. 
\smallskip
\item{[Ag]} Agmon S., {\sl Lectures on Elliptic Boundary Value Problems}, Van Nostrand, Princeton, N.J., 1965.\smallskip
\item{[AT]} Adachi S., Tanaka K., 
{\sl Trudinger type inequalities in $\R^N$ and their best exponents}, 
 Proc. Amer. Math. Soc. {\bf 128} (2000), 2051-2057.\smallskip 
\item{[BGS]} Banica V., Gonz\'alez M.d.M., S\'aez M., {\sl Some constructions for the fractional Laplacian on noncompact manifolds}, Rev. Mat. Iberoam. {\bf 31} (2015), 681-712.\smallskip
\item{[BM]} Battaglia L., Mancini G., {\sl Remarks on the Moser-Trudinger inequality},   Adv. Nonlinear Anal. {\bf 2} (2013), 389-425.\smallskip 
\item{[BS]} Bennett C., Sharpley R., {\sl Interpolation of Operators}, Pure and  Applied  Math., 129, Academic Press, Inc., Boston, MA, 1988.\smallskip
\item{[C]} Cao D. M., 
{\sl Nontrivial solution of semilinear elliptic equation with critical exponent in $\R^2$}, 
Comm. Partial Differential Equations {\bf  17} (1992),  407-435. \smallskip

\item{[CK]} Cohl H.S., Kalnins E.G.,  {\sl Fourier and Gegenbauer expansions for a fundamental solution of the Laplacian in the hyperboloid model of hyperbolic geometry}, J. Phys. A: Math. Theor. {\bf 45} (2012), 1-33.\smallskip

\item{[do\'O]}  do \'O J.M.B. {\sl $N$-Laplacian equations in $\R^N$ with critical growth}, Abstr. Appl. Anal. {\bf 2} (1997), 301-315.\smallskip

\item{[FM1]} Fontana L., Morpurgo C., {\sl Adams inequalities on measure spaces}, Adv. Math. {\bf 226} (2011), 5066-5119.
\smallskip
\item{[FM2]} Fontana L., Morpurgo C., {\sl Sharp exponential integrability for critical Riesz potentials and fractional Laplacians on $\Rn$}, Nonlinear Anal. {\bf 167} (2018), 85-122.\smallskip
\item{[FM3]} Fontana L., Morpurgo C., {\sl Sharp Adams and Moser-Trudinger inequalities on $\Rn$ and other spaces of infinite measure}, preprint (2015),	arXiv:1504.04678. \smallskip
\item{[LL]} Lam N., Lu G., {\sl A new approach to sharp Moser-Trudinger and Adams type inequalities: a rearrangement-free argument}, 
J. Differential Equations {\bf 255} (2013), 298-325.\smallskip 
\item{[LR]} Li Y., Ruf B., {\sl A sharp Trudinger-Moser type inequality for unbounded domains in $R^n$},   
 Indiana Univ. Math. J. {\bf 57} (2008),  451-480.\smallskip
 \item{[Lo]} Lorenzi A., {\sl On elliptic equations with piecewise constant coefficients. II}, Ann. Scuola Norm. Sup. Pisa Cl. Sci. {\bf 26} (1972), 839-870.\smallskip 
\item{[LT1]}  Lu G., Tang  H., {\sl  Best constants for Moser-Trudinger inequalities on high dimensional hyperbolic spaces}, Adv. Nonlinear Stud. {\bf 13} (2013),  1035-1052. \smallskip
\item{[LT2]}  Lu G., Tang  H., {\sl  Sharp Moser-Trudinger Inequalities on Hyperbolic
Spaces with Exact Growth Condition}, J. Geom. Anal. {\bf 26} (2016), 837-857.\smallskip

 \item{[MS]} Mancini G., Sandeep K., {\sl  Moser-Trudinger inequality on conformal discs}, Commun. Contemp. Math. {\bf 12} (2010),  1055-1068. \smallskip
\item{[MST]} Mancini G., Sandeep K., Tintarev C.,
{\sl Trudinger-Moser inequality in the hyperbolic space $\H^N$}, 
 Adv. Nonlinear Anal. {\bf 2} (2013),  309-324.\smallskip

\item{[MaS]} Masmoudi N., Sani F., {\sl Higher order Adams' inequality with the exact growth condition}, Commun. Contemp. Math. {\bf20} (2018), 33 pp. 
\smallskip
\item{[OW]} Ortner N., Wagner P., {\sl Fundamental solutions of linear partial differential operators}, Springer 2015.\smallskip
\item{[P]} Panda R. {\sl 
Nontrivial solution of a quasilinear elliptic equation with critical growth in $\R^n$}, Proc. Indian Acad. Sci. Math. Sci. {\bf 105} (1995), 425-444.\smallskip 
\item{[QS]} Quittner P., Souplet P., {\sl
Superlinear parabolic problems.
Blow-up, global existence and steady states. } Birkhäuser Advanced Texts: Basler Lehrbücher. Birkhäuser Verlag, Basel, 2007.
\smallskip
\item{[RiSt]} Ricci F., Stein E.M.,  {\sl  Harmonic analysis on nilpotent groups  and singular integrals I. Oscillatory integrals}, J. Funct. Anal. {\bf 73} (1987), 179-194.\smallskip
\item{[R]} Ruf B., {\sl A sharp Trudinger--Moser type inequality for unbounded domains in $\R^2$}, J. Funct. Anal. {\bf 219} (2005), 340-367.\smallskip

\item{[RS]} Ruf B., Sani F., {\sl
 Sharp Adams-type inequalities in $\R^n$},  
Trans. Amer. Math. Soc. {\bf 365} (2013), 645-670.\smallskip 
\item{[S]} Stein E.M., {\sl
Singular integrals and differentiability properties of functions}, 
Princeton Mathematical Series {\bf 30}, Princeton University Press, 1970.\smallskip
\item{[So]} Souplet P., {\it Geometry of unbounded domains, Poincaré inequalities and stability in semilinear parabolic equations}, 
Comm. Partial Differential Equations {\bf 24} (1999),  951-973.
\item{[Tar]} Tarsi C., 
{\sl Adams' inequality and limiting Sobolev embeddings into Zygmund spaces}, 
Potential Anal. {\bf 37} (2012),  353-385.\smallskip 
\item{[Tat]}  Tataru D., {\sl Strichartz estimates in the hyperbolic space and global existence for the semilinear wave equation}, Trans. Amer. Math. Soc. {\bf 353} (2001), 795-807.\smallskip
\item{[Tr]} Trudinger N.S., {\sl 
On imbeddings into Orlicz spaces and some applications}, 
J. Math. Mech. {\bf17} (1967), 473-483. 
\smallskip
\item{[Y]} Yang Y., {\sl 
Trudinger-Moser inequalities on complete noncompact Riemannian manifolds},  
J. Funct. Anal. {\bf 263} (2012),  1894-1938.\smallskip

\smallskip\smallskip
\noin Luigi Fontana \hskip19em Carlo Morpurgo

\noin Dipartimento di Matematica ed Applicazioni \hskip5.5em Department of Mathematics 
 
\noin Universit\'a di Milano-Bicocca\hskip 12.8em University of Missouri

\noin Via Cozzi, 53 \hskip 19.3em Columbia, Missouri 65211

\noin 20125 Milano - Italy\hskip 16.6em USA 
\smallskip\noin luigi.fontana@unimib.it\hskip 15.3em morpurgoc@missouri.edu

\end